\documentclass[12pt]{article}
\usepackage[intlimits]{amsmath}
\usepackage{amssymb}
\usepackage{mathrsfs}
\usepackage[russian]{babel}
\usepackage[dvips]{graphicx}
\usepackage{longtable}
\usepackage{indentfirst}

\title{Дистанционные графы, имеющие большое хроматическое число и не содержащие клик или циклов заданного размера\footnote{Работа выполнена при
финансовой поддержке гранта РФФИ 09-01-00294, гранта Президента РФ  
МД-8390.2010.1, гранта поддержки ведущих научных школ НШ-8784.2010.1.}}

\author{Е.Е. Демёхин, А.М. Райгородский, О.И. Рубанов}

\date{}

\begin{document}

\maketitle

\begin{flushright}

УДК 519.174.7+519.176+519.175.4

\end{flushright}

\section{Введение и постановки основных задач}

Настоящая работа посвящена решению ряда экстремальных задач для дистанционных графов в евклидовом пространстве 
$ {\mathbb R}^n $ растущей размерности. Напомним, что 
{\it дистанционный граф} (или {\it граф расстояний}) -- это такой граф, у которого вершины являются элементами некоторого метрического 
пространства $ X $ с метрикой $ \rho $, а ребра -- те или иные пары вершин, отстоящих друг от друга на расстояние из данного множества 
$ {\cal A} \subseteq {\mathbb R}_{+} $. Иными словами, граф $ G = (V,E) $ дистанционный, если 
$$
V \subseteq X, ~~~ E \subseteq \{\{x,y\}: ~ x,y \in V, ~ \rho(x,y) \in {\cal A}\}. 
$$

Дистанционные графы играют огромную роль в комбинаторной геометрии. Например, они естественным образом возникают в связи с классической 
задачей о {\it хроматическом числе пространства} $ {\mathbb R}^n $. Эта задача состоит в отыскании наименьшего количества $ \chi({\mathbb R}^n) $
цветов, в которые можно так покрасить все точки $ {\mathbb R}^n $, чтобы между одноцветными точками не было расстояния 1: 
$$
\chi({\mathbb R}^n) = \min \left\{\chi:~ {\mathbb R}^n = V_1 \bigsqcup \ldots \bigsqcup V_{\chi}, ~ \forall ~i ~\forall ~ {\bf x}, {\bf y} \in V_i ~~
|{\bf x}-{\bf y}| \neq 1\right\}.
$$
Здесь $ |{\bf x}-{\bf y}| $ -- евклидово расстояние между векторами. Понятно, причем здесь и дистанционные графы. Действительно, возьмем граф 
$ \mathfrak{G} = (\mathfrak{V}, \mathfrak{E}) $, у которого 
$$
\mathfrak{V} = {\mathbb R}^n, ~~~ \mathfrak{E} = \{\{{\bf x},{\bf y}\}:~ |{\bf x}-{\bf y}| = 1\}.
$$
Рассмотрим его {\it хроматическое число} $ \chi(\mathfrak{G}) $ (минимальное количество цветов, в которые можно так покрасить все вершины графа, чтобы у него 
не было ребер с одноцветными концами). Ясно, что $ \chi({\mathbb R}^n) = \chi(\mathfrak{G}) $. Более того, по теореме П. Эрдеша -- Н.Г. де Брёйна
(см. \cite{ErBr}), если величина $ \chi({\mathbb R}^n) $ конечна, то она достигается на некотором конечном подграфе графа $ \mathfrak{G} $. 
Однако совсем легко видеть, что, скажем, $ \chi({\mathbb R}^n) \le (n+1)^n $, и, стало быть, с точки зрения хроматического числа пространства, 
достаточно ограничиться изучением конечных графов расстояний в $ {\mathbb R}^n $. 

Перечислим основные результаты о хроматических числах дистанционных графов в различных размерностях. Нас будет интересовать в первую очередь 
максимальное значение хроматического числа, т.е. величина $ \chi({\mathbb R}^n) $. 

\begin{itemize}  

\item{$ \chi({\mathbb R}^1) = 2 $. Этот результат очевиден.}

\item{$ 4 \le \chi({\mathbb R}^2) \le 7 $. Оба результата довольно простые. Доказательства их можно найти в работах 
\cite{MM}, \cite{Had} и в брошюре \cite{Rai1}.}

\item{$ 6 \le \chi({\mathbb R}^3) \le 15 $. Нижняя оценка принадлежит О. Нечуштану (см. \cite{Nech}), а верхняя -- Д. Кулсону (см. \cite{Coul}).}

\item{$ 7 \le \chi({\mathbb R}^4) \le 54 $. Нижняя оценка принадлежит К. Кантвеллу (см. \cite{Cant}) и 
Л.Л. Иванову (см. \cite{Iva}), а верхняя -- Г. Тоту и Р. Радойчичу (см. \cite{TR}).}

\item{$ (1.239...+o(1))^n \le \chi({\mathbb R}^n) \le (3+o(1))^n $. Нижняя оценка принадлежит А.М. Райгородскому (см. \cite{Rai2}), а верхняя -- 
Д. Ларману и К.А. Роджерсу (см. \cite{LR}). Заметим, что приведенная нижняя оценка является модификацией первоначальной экспоненциальной 
оценки П. Франкла и Р.М. Уилсона, имевшей вид $ \chi({\mathbb R}^n) \ge (1.207...+o(1))^n $ (см. \cite{FW}).}

\end{itemize}  

Таким образом, мы видим, что с ростом размерности хроматическое число пространства (максимальное значение хроматического числа графа расстояний 
в $ {\mathbb R}^n $) растет экспоненциально. 

Отметим также, что для дистанционных графов, которые мы сейчас рассматриваем, множество $ {\cal A} $, порождающее ребра, имеет тривиальный вид 
$ {\cal A} = \{1\} $. Ясно, что в рамках задачи о хроматических числах мы с тем же успехом можем взять любое $ {\cal A} = \{a\} $: от этого
ничего не изменится. Заметим, наконец, что, поскольку мы {\it запрещаем} точкам одного цвета отстоять друг от друга на расстояние 1 (или $ a $), 
величину $ a $ зачастую называют {\it запрещенным расстоянием}. 

Перейдем к постановкам экстремальных задач, о которых пойдет речь в настоящей работе. В 1959 году Эрдеш (см. \cite{Erd1}) доказал красивый и 
довольно неожиданный факт. 

\vskip+0.2cm

\noindent {\bf Теорема 1.} {\it Пусть $ k \ge 2 $, $ l \ge 2 $ -- произвольные натуральные числа. Тогда существует граф, у которого 
хроматическое число больше $ k $, а длина минимального простого цикла больше $ l $.}

\vskip+0.2cm

Иными словами, бывают графы $ G $ со сколь угодно большим хроматическим числом $ \chi(G) $ и {\it обхватом} $ g(G) $ (длиной кратчайшего цикла). 

В 1976 году Эрдеш, отталкиваясь от своей же теоремы 1, сформулировал такой вопрос (см. \cite{Erd2}): {\it существует ли 
дистанционный граф на плоскости, имеющий хроматическое число 4 и не содержащий треугольников?} Отметим, что конструкция из \cite{MM}, 
дающая оценку $ \chi({\mathbb R}^2) \ge 4 $, содержит четыре треугольника, и это по существу. На вопрос Эрдеша был получен 
утвердительный ответ. Уже в 1979 году Н. Уормалд построил плоский граф расстояний с хроматическим числом 4 и без треугольников (см. \cite{Worm}).
Правда, у этого графа было 6448 вершин. В 1996 году П. О'Доннелл и Р. Хохберг предложили аналогичный граф всего с двадцатью тремя вершинами (см. 
\cite{Hoch}). Наконец, в 2000 году О'Доннелл доказал существование графа расстояний на плоскости с хроматическим числом 4 и наперед заданным 
обхватом (см. \cite{Don1}, \cite{Don2}). 

Нас интересуют обобщения задачи Эрдеша на случай размерности $ n \ge 3 $. Разумеется, здесь речь идет о построении графов с хроматическим числом, 
максимально близким к наибольшему из известных. Так, в размерности 3 имеет смысл говорить о хроматическом числе 5 или 6. В растущей размерности
правильным аналогом четверки служит экспонента. Остается лишь ввести основные величины, которые нам предстоит изучать. 

Итак, пусть $ \omega(G) $ -- это число вершин в максимальном полном подграфе графа $ G = (V,E) $:
$$
\omega(G) = \max \{|W|:~ W \subseteq V, ~ \forall~ x,y \in W ~~ \{x,y\} \in E\}.
$$
Поскольку любой полный подграф называется {\it кликой}, величину $ \omega(G) $ называют {\it кликовым числом} графа. Утверждение "$ \omega(G) < 3 $", 
очевидно, равносильно утверждению об отсутствии треугольников в $ G $. Положим 
$$
\zeta_{{\rm clique}}(k) = \sup \left\{\zeta: ~ \exists ~ {\rm функция} ~ \delta = \delta(n), ~ {\rm такая}, ~ {\rm что} ~ \lim\limits_{n \to \infty} \delta(n) = 0 ~
\right.
$$
$$
\left. \phantom{\lim\limits_{n \to \infty}} 
{\rm и} ~ \forall ~ n ~ \exists ~ G ~ {\rm в} ~ {\mathbb R}^n, ~ {\rm у} ~ {\rm которого} ~ \omega(G) < k, ~ \chi(G) \ge (\zeta+\delta(n))^n\right\},
$$
$$
\zeta_{{\rm girth}}(k) = \sup \left\{\zeta: ~ \exists ~ {\rm функция} ~ \delta = \delta(n), ~ {\rm такая}, ~ {\rm что} ~ \lim\limits_{n \to \infty} \delta(n) = 0 ~
\right.
$$
$$
\left. \phantom{\lim\limits_{n \to \infty}} 
{\rm и} ~ \forall ~ n ~ \exists ~ G ~ {\rm в} ~ {\mathbb R}^n, ~ {\rm у} ~ {\rm которого} ~ g(G) > k, ~ 
\chi(G) \ge (\zeta+\delta(n))^n\right\}.
$$
Разумеется, для некоторых $ k $ эти величины могут оказаться меньше единицы, и тогда смысла в их рассмотрении нет. Основная цель настоящей 
работы состоит в том, чтобы доказать неравенства $ \zeta_{{\rm clique}}(k) > 1 $ для всех $ k $ и, более того, найти оптимальные значения в 
правых частях этих неравенств. Аналогичное исследование мы проведем и для случая величины $ \zeta_{{\rm girth}}(k) $. 

Дальнейшая структура статьи будет следующей. В разделе 2 мы подробно изучим величину $ \zeta_{{\rm clique}}(3) $. В разделе 3 мы получим 
серию оценок для $ \zeta_{{\rm clique}}(k) $ с произвольным $ k $, исходя из некоторого общего подхода, связанного с теорией кодирования. 
В разделе 4 мы также получим оценки для $ \zeta_{{\rm clique}}(k) $ с произвольным $ k $, но при этом используем принципиально другой 
подход, нежели в разделе 3, -- подход вероятностный. В том же разделе 4 мы аккуратно сопоставим все полученные результаты. Раздел 5 мы 
посвятим оценкам величины $ \zeta_{{\rm girth}}(k) $. Наконец, в шестом разделе мы обобщим полученные результаты на случай нескольких 
запрещенных расстояний (дистанционных графов в $ {\mathbb R}^n $ с множествами запретов $ {\cal A} $, $ |{\cal A}| > 1 $). 

Завершая раздел, заметим, что историю задачи о хроматическом числе пространства и различные результаты относительно дистанционных 
графов можно найти в обзорах \cite{Rai3}, \cite{Szek}, \cite{Soi1} и книгах \cite{Rai1}, \cite{BMP}, \cite{PA}, \cite{KW}, \cite{Soi2}, 
\cite{Rai4}.

\section{Оценка величины $ \zeta_{{\rm clique}}(3) $}

Этот раздел мы построим следующим образом. В параграфе 2.1 мы сформулируем наш основной результат. В параграфе 2.2 мы его докажем. 
Параграфы 2.3 и 2.4 будут посвящены обсуждению возможности уточнения основного результата. 

\subsection{Формулировка основного результата}

Сперва заметим, что еще в 2007 году О.И. Рубанов доказал следующую теорему (см. \cite{Rub}). 

\vskip+0.2cm

\noindent {\bf Теорема 2.} {\it Существует такой граф расстояний $ G $ в 
$ {\mathbb R}^3 $, что $ \chi(G) \ge 5 $ и $ \omega(G) \le 3 $.}

\vskip+0.2cm

В этом году Рубанову удалось заменить оценку $ \omega(G) \le 3 $ на $ \omega(G) \le 2 $. Данный результат пока не опубликован. Что же 
до основного нашего вопроса -- вопроса о графах без треугольников в $ {\mathbb R}^n $, $ n \to \infty $, и их хроматических числах, -- 
то принципиально новый результат по нему содержится в теореме 3. 

\vskip+0.2cm

\noindent {\bf Теорема 3.} {\it Имеет место оценка 
$$ 
\zeta_{{\rm clique}}(3) \ge 2 \cdot \left(\frac{1}{3}\right)^{1/3} \cdot \left(\frac{2}{3}\right)^{2/3} = 1.058... 
$$}

\vskip+0.2cm

Теорема 3 означает, что в растущей размерности есть графы расстояний без треугольников, у которых, тем не менее, экспоненциально 
большие хроматические числа. Конечно, величина $ 1.058 $ меньше и величины $ 1.239 $, и величины $ 1.207 $ (см. раздел 1), но она все-таки 
строго больше единицы. Более того, из теоремы 3 моментально следует осмысленность изучения {\it всех} величин $ \zeta_{{\rm clique}}(k) $, 
т.к. мы теперь знаем, что каждая из этих величин никак не меньше, чем $ 1.058 $. 

В следующем параграфе мы докажем теорему 3. 

\subsection{Доказательство теоремы 3}

Зафиксируем размерность $ n $. Положим $ k = [0.5n] $. Рассмотрим минимальное простое число $ p $, такое, что $ k - p < \frac{n}{6}-1 $. 
Известно, что с ростом $ n $ справедлива асимптотика $ p = p(n) \sim \frac{n}{3} $. Это обусловлено тем, что между $ x $ и 
$ x + O\left(x^{0.525}\right) $ всегда есть простое число (см. \cite{BHP}).

Рассмотрим граф $ G = (V,E) $, у которого 
$$
V = \{{\bf x} = (x_1, \dots, x_n): ~ \forall ~i~~ x_i \in \{0,1\}, ~ x_1 + \ldots + x_n = k\}, ~~
E = \{\{{\bf x},{\bf y}\}:~ {\bf x}, {\bf y} \in V, ~ |{\bf x}-{\bf y}| = \sqrt{2p}\}.
$$
Иными словами, здесь запрещенное расстояние -- это величина $ a = \sqrt{2p} $, которая однозначно соответствует скалярному 
произведению $ ({\bf x},{\bf y}) = k-p $. 

Прежде всего заметим, что в графе $ G $ нет треугольников. Действительно, предположим, что треугольник с вершинами $ {\bf x}, {\bf y}, {\bf z} $ 
нашелся. Обозначим через $ {\cal R}_n $ множество $ \{1, \dots, n\} $ номеров координат каждого из векторов $ {\bf t} \in V $. Пусть 
$ A_{{\bf x}}, A_{{\bf y}}, A_{{\bf z}} \subset {\cal R}_n $ -- множества позиций, на которых расположены единицы векторов 
$ {\bf x}, {\bf y}, {\bf z} $ соответственно. Поскольку пара $ {\bf x}, {\bf y} $ образует ребро, $ |A_{{\bf x}} \cup A_{{\bf y}}| = 
k+p $. На позициях из $ A_{{\bf x}} \cup A_{{\bf y}} $ вектор $ {\bf z} $ не может иметь больше, чем $ 2k-2p $, единиц. Значит, 
остальные $ k - (2k-2p) = 2p-k $ его единиц должны разместиться на позициях из $ {\cal R}_n \setminus (A_{{\bf x}} \cup A_{{\bf y}}) $. 
Количество этих позиций равно $ n - k-p $. Если мы покажем, что $ n-k-p < 2p-k $, то мы придем к противоречию с предположением о наличии 
треугольника в $ G $. Итак, 
$$
(n-k-p)-(2p-k) = n - 3p < n - 3 \left(k-\frac{n}{6}+1\right) \le n - 3 \left(\frac{n}{2} - \frac{n}{6}\right) = 0.
$$
Искомое неравенство установлено, и, значит, в $ G $ на самом деле отсутствуют треугольники.  

Перейдем к оценке хроматического числа. Хорошо известно неравенство $ \chi(G) \ge \frac{|V|}{\alpha(G)} $, где $ \alpha(G) $ -- это 
величина, в некотором смысле двойственная величине $ \omega(G) $ (отсюда и обозначение):
$$
\alpha(G) = \max \{|W|:~ W \subseteq V, ~ \forall~ x,y \in W ~~ \{x,y\} \not \in E\}.
$$
Эта величина называется {\it числом независимости} графа; в свою очередь, все множества вершин графа, свободные от его ребер, называются 
{\it независимыми}. Таким образом, число независимости -- это мощность самого большого независимого множества вершин. 

В нашем случае, очевидно, $ |V| = C_n^k $. Значит, нам нужно как можно точнее оценить сверху $ \alpha(G) $. Это делается с помощью 
линейно-алгебраического метода в комбинаторике (см. \cite{Rai4}, \cite{BF}). А именно, каждому вектору $ {\bf x} \in V $ 
мы сопоставляем некоторый полином $ P_{{\bf x}} \in {\mathbb Z}/p {\mathbb Z}[y_1, \dots, y_n] $. В данном случае мы полагаем 
$$
P_{{\bf x}}({\bf y}) = \prod_{i=0, i \not \equiv k \pmod p}^{p-1} (i-({\bf x},{\bf y})), ~~~ {\bf y} = (y_1, \dots, y_n). 
$$
Иными словами, произведение берется по всем наименьшим неотрицательным вычетам по модулю $ p $, кроме остатка от деления $ k $ на $ p $. 

Полиномы подобраны так, чтобы выполнялось следующее 

\vskip+0.2cm

\noindent {\bf Свойство 1.} {\it Для любых $ {\bf x}, {\bf y} \in V $ условие $ ({\bf x},{\bf y}) \equiv k \pmod p $ равносильно 
условию $ P_{{\bf x}}({\bf y}) \not \equiv 0 \pmod p $.}

\vskip+0.2cm

Свойство очевидно, и мы его не доказываем. Ясно, далее, что степень каждого полинома не превосходит $ p-1 $. Сейчас мы еще упростим 
полиномы, причем так, что свойство 1 для новых полиномов останется верным. Раскроем в данном полиноме $ P_{{\bf x}} $ скобки. Получится 
линейная комбинация мономов, каждый из которых имеет вид
$$
y_{i_1}^{\alpha_1} \cdot \ldots \cdot y_{i_q}^{\alpha_q}, ~~~ q \le p - 1. 
$$
Формально заменим всякий такой моном мономом 
$$
y_{i_1} \cdot \ldots \cdot y_{i_q}.
$$
Поскольку в свойстве 1 фигурируют лишь векторы $ {\bf y} \in V $, т.е. векторы, все координаты которых суть 0 или 1, то значение 
полинома не поменяется, а стало быть, свойство 1 для новых полиномов $ P_{{\bf x}}' $ снова выполнено. 

Полиномы $ P_{{\bf x}}' $ хороши тем, что, хотя их степени по-прежнему могут достигать значения $ p-1 $, все они содержатся в 
пространстве полиномов размерности $ \sum\limits_{i=0}^{p-1} C_n^i < n C_n^p $. На сравнительной малости этой размерности мы и сыграем.

Рассмотрим произвольное независимое множество $ W = \{{\bf x}_1, \dots, {\bf x}_s\} \subset V $. Иначе говоря, для любых $ i,j $
выполнено $ ({\bf x}_i,{\bf x}_j) \neq k-p $. Возьмем полиномы $ P_{{\bf x}_1}', \dots, P_{{\bf x}_s}' $. Предположим, что для 
некоторых $ c_1, \dots, c_s \in {\mathbb Z}/p {\mathbb Z} $ имеет место сравнение
$$
c_1 P_{{\bf x}_1}'({\bf y}) + \ldots + c_s P_{{\bf x}_s}'({\bf y}) \equiv 0 \pmod p ~~~ \forall ~ {\bf y} \in V.
$$
Пусть $ {\bf y} = {\bf x}_i $, $ i $ -- любое. Тогда 
$$ 
({\bf x}_i,{\bf y}) = ({\bf x}_i,{\bf x}_i) = k \equiv k \pmod p.
$$
Значит, ввиду свойства 1, $ P_{{\bf x}_i}'({\bf y}) \not \equiv 0 \pmod p $. В то же время $ ({\bf x}_j,{\bf y}) \neq k-p $. Разумеется, 
$ ({\bf x}_j,{\bf y}) < k $. Но очевидно, что $ k - 2p < 0 $, т.е. $ ({\bf x}_j,{\bf y}) > k-2p $. В итоге $ ({\bf x}_j,{\bf y}) \not \equiv k
\pmod p $, а стало быть, опять же за счет свойства 1, $ P_{{\bf x}_j}'({\bf y}) \equiv 0 \pmod p $ для всех $ j \neq i $. Поскольку $ p $ -- 
простое, имеем $ c_i \equiv 0 \pmod p $ при любом $ i $. 

Приведенное выше рассуждение показывает, что полиномы $ P_{{\bf x}_1}', \dots, P_{{\bf x}_s}' $ линейно независимы над 
$ {\mathbb Z}/p {\mathbb Z} $. Следовательно, $ s $ не превосходит размерности пространства, порожденного этими полиномами, а она, как мы помним, 
не больше, чем $ n C_n^p $. Множество $ W $ было произвольным, так что выполнено и неравенство $ \alpha(G) \le n C_n^p $. 

Итак, 
$$
\chi(G) \ge \frac{C_n^k}{n C_n^p}.
$$
Пользуемся тем фактом, что $ p \sim \frac{n}{3} $, и стандартными выкладками с участием формулы Стирлинга. Получаем 
$$
C_n^k = (2+o(1))^n, ~~~ nC_n^p = \left(\frac{1}{(1/3)^{1/3} (2/3)^{2/3}}+o(1)\right)^n, ~~~ 
\zeta_{{\rm clique}}(3) \ge 2 \cdot \left(\frac{1}{3}\right)^{1/3} \cdot \left(\frac{2}{3}\right)^{2/3} = 1.058...
$$

Теорема доказана. 

\subsection{Общая конструкция с (0,1)-векторами}

В доказательстве теоремы 3 мы использовали граф, вершины которого представляли собой (0,1)-векторы с фиксированным количеством $ k $ единиц. 
При этом $ k $ было выбрано равным $ [0.5n] $. При внимательном изучении всего рассуждения, приведенного в \S 2.2, становится ясно, что с 
самого начала можно было брать следующую более общую конструкцию. 

Пусть $ k \in {\mathbb N} $, $ k \le n $. Пусть также $ p $ -- любое простое число. Положим 
$$
V = \{{\bf x} = (x_1, \dots, x_n): ~ \forall ~i~~ x_i \in \{0,1\}, ~ x_1 + \ldots + x_n = k\}, ~~
E = \{\{{\bf x},{\bf y}\}:~ {\bf x}, {\bf y} \in V, ~ |{\bf x}-{\bf y}| = \sqrt{2p}\}.
$$
При $ p < \frac{k}{2} $ треугольники в конструкции точно есть. При $ p \ge \frac{k}{2} $
условием отсутствия треугольников в конструкции по-прежнему служит неравенство 
$$
n - k - p < 2p - k.
$$
Иными словами, $ p > \frac{n}{3} $ при любом $ k $. Условием применимости линейно-алгебраического метода является, в свою очередь, 
неравенство 
$$ 
k - 2p < \min_{{\bf x},{\bf y}} ({\bf x},{\bf y}) = \max \{0, 2k-n\}.
$$
Если $ k \le \frac{2n}{3} $ и $ p > \frac{n}{3} $, то 
$$
k-2p < \frac{2n}{3} - \frac{2n}{3} = 0 \le \max \{0, 2k-n\}.
$$
Если же $ \frac{2n}{3} < k \le n $ и $ p > \frac{n}{3} $, то 
$$
k-2p < n - \frac{2n}{3} = \frac{n}{3} < \max \{0, 2k-n\} = 2k-n.
$$
Таким образом, условие отсутствия треугольников всегда влечет применимость линейной алгебры.

Получается оценка 
$$
\chi(G) \ge \frac{C_n^k}{n C_n^p}, 
$$
которая, очевидно, тем больше, чем меньше $ p $. Следовательно, мы возвращаемся к ситуации $ p \sim \frac{n}{3} $. Более того, 
уже при данном $ p $ оценка тем лучше, чем большее значение принимает величина $ C_n^k $. Так мы возвращаемся к $ k = [0.5n] $. 

Мы показали, что параметры в теореме 3 были выбраны оптимально с точки зрения предложенного метода.

\subsection{Дальнейшие возможные обобщения}

Графы, множества вершин которых состоят из (0,1)-векторов, в науке о хроматических числах впервые были использованы в работе \cite{FW}, 
где, как мы уже говорили в разделе 1, была доказана оценка $ \chi({\mathbb R}^n) \ge (1.207...+o(1))^n $. Впоследствии рассматривались 
и графы более общего вида. Например, наилучшая известная оценка хроматического числа пространства (с констатой 1.239 в основании 
экспоненты) построена на рассмотрении графов с вершинами в (-1,0,1)-векторах (см. \cite{Rai2}). Изучались и графы, вершины 
которых являются произвольными элементами решетки $ {\mathbb Z}^n $ (см. \cite{Rai5}, \cite{Rai6}, \cite{Sh1}, \cite{Sh2}). 

В этой связи уместно было бы попробовать строить и графы без треугольников хотя бы на (-1,0,1)-векторах. К сожалению, к успеху 
это не приводит. Нам удалось доказать 

\vskip+0.2cm 	

\noindent {\bf Предложение 1.} {\it Пусть $ n \in {\mathbb N} $. Пусть, далее, 
$ k_1, k_{-1}, k_0 $ -- произвольные натуральные числа, в сумме дающие $ n $. Пусть, наконец, $ a $ -- произвольное число. Положим 
$$
V = \{{\bf x} = (x_1, \dots, x_n): \forall ~i~ x_i \in \{-1,0,1\}, ~ |\{i:x_i = 1\}| = k_1, ~ |\{i:x_i = -1\}| = k_{-1}, ~
|\{i:x_i = 0\}| = k_0\}, 
$$
$$
E = \{\{{\bf x},{\bf y}\}:~ {\bf x}, {\bf y} \in V, ~ |{\bf x}-{\bf y}| = a\}.
$$
Тогда при любых значениях параметров $ k_1, k_{-1}, k_0 $ и $ a $, которые гарантируют отсутствие треугольников в конструкции, 
нижняя оценка 
хроматического числа графа $ G = (V,E) $, полученная с помощью линейно-алгебраического метода,  
не превосходит величины $ (c+o(1))^n $, где 
$$ 
c \le 2 \cdot \left(\frac{1}{3}\right)^{1/3} \cdot \left(\frac{2}{3}\right)^{2/3} = 1.058... 
$$}

\vskip+0.2cm

Доказательство предложения 1 основано на подробном разборе ряда случаев (соотношений между параметрами), и, поскольку оно, с одной 
стороны, громоздко, с другой стороны, просто, а главное, не позволяет получить новый результат в исходной задаче, мы его здесь 
не приводим.  

Для более сложносоставленных векторов целочисленной решетки довольно тяжело выписывать условия отсутствия треугольников. По-видимому, 
исчерпывающий результат в этих случаях возможен лишь за счет применения численных методов. В настоящей работе мы этим не занимались.

\section{Оценка величины $ \zeta_{{\rm clique}}(k) $ с помощью границ для равновесных кодов}

Этот раздел мы разобьем на три параграфа. В первом параграфе мы обсудим общую идею подхода и сформулируем наш результат. Во втором 
параграфе мы докажем этот результат. В третьем параграфе мы сделаем несколько замечаний и приведем таблицу оценок для разных 
значений $ k $. 

\subsection{Идея подхода и формулировка результата}

Мы снова хотим использовать графы, у которых вершины -- (0,1)-векторы с фиксированным числом единиц. В теории кодирования совокупности 
таких векторов называются {\it равновесными кодами}, и про равновесные коды известно достаточно много. Например, доказан ряд верхних 
оценок для числа кодовых слов (т.е. как раз (0,1)-векторов) в равновесных кодах с данной нижней границей для расстояния Хэмминга. В нашем 
случае нижняя граница для расстояния Хэмминга -- это, по сути, величина запрещенного расстояния. Поэтому ясно, что применение известных 
результатов для кодов может сильно помочь в нашей задаче. И действительно, удается доказать следующую общую теорему. 

\vskip+0.2cm

\noindent {\bf Теорема 4.} {\it Пусть дано натуральное число $ k \ge 3 $. Пусть, далее, $ a $ и $ b $ -- произвольные вещественные числа, 
удовлетворяющие ограничениям
$$
a \in (0,1), ~~ b \in (0,1), ~~ b > 2a, ~~ b^2 > a, ~~ \left \lceil \frac{b-a}{b^2-a} \right \rceil < k.
$$
Тогда 
$$
\zeta_{{\rm clique}}(k) \ge \frac{(b-a)^{b-a} (1-b+a)^{1-b+a}}{b^b (1-b)^{1-b}}.
$$} 

\vskip+0.2cm

Результаты для конкретных значений $ k $, полученные путем оптимизации по $ a $ и $ b $, мы приведем в параграфе 3.3. 
Доказательство -- в следующем параграфе. 

\subsection{Доказательство теоремы 4}

Положим $ l = [bn] $. Возьмем минимальное простое число $ p $, для которого $ l-p < [an] $ и 
$$
\frac{pn}{l^2-ln+pn} \le \frac{b-a}{b^2-a}.
$$
Тогда $ p \sim (b-a)n $. 
При этом за счет неравенства $ b > 2a $ мы имеем также оценку $ l-2p < 0 $. Рассмотрим граф $ G=(V,E) $, у которого 
$$
V = \{{\bf x} = (x_1, \dots, x_n): ~ \forall ~i~~ x_i \in \{0,1\}, ~ x_1 + \ldots + x_n = l\}, ~~
E = \{\{{\bf x},{\bf y}\}:~ {\bf x}, {\bf y} \in V, ~ |{\bf x}-{\bf y}| = \sqrt{2p}\}.
$$

Справедлива следующая теорема.  

\vskip+0.2cm

\noindent {\bf Теорема 5.} {\it Максимальное число двоичных $ n $-мерных векторов веса 
$ l $ (т.е. с $ l $ единицами), отстоящих друг от друга на расстояние Хэмминга, не меньшее $ 2 \delta $, 
не превосходит величины 
$$
\left \lceil \frac{\delta n}{l^2 - ln + \delta n} \right \rceil, 
$$ 
коль скоро знаменатель указанной дроби больше нуля.}  

\vskip+0.2cm
 
Доказательство теоремы можно найти в книге \cite{WS}. Эта теорема позволяет утверждать, что $ \omega(G) < k $. 
Действительно, пусть $ H \subset G $ -- некоторая клика. Тогда ее вершины отстоят друг от друга на расстояние 
Хэмминга, равное $ 2p $. Значит, $ \delta = p $ и размер клики $ H $ никак не больше, чем
$$
\left \lceil \frac{pn}{l^2-ln+pn} \right \rceil \le \left \lceil \frac{b-a}{b^2-a} \right \rceil < k,
$$
ведь по условию $ b^2 > a $. 

Далее, линейно-алгебраический метод (применимый ввиду неравенства $ l - 2p < 0 $) дает оценку 
$$
\chi(G) \ge \frac{C_n^l}{n C_n^p}.
$$
Поскольку $ l \sim bn $, а $ p \sim (b-a)n $, получаем за счет стандартных выкладок, что 
$$
\chi(G) \ge \left(\frac{(b-a)^{b-a} (1-b+a)^{1-b+a}}{b^b (1-b)^{1-b}} + o(1)\right)^n.
$$

Теорема доказана.

\subsection{Комментарии и таблица оценок}

Прежде всего отметим, что при $ k \le 4 $ результат теоремы 4 слабее результата теоремы 3. Это связано с тем, что оценка из 
теории кодирования универсальна и мы, по сути, слегка уточнили ее при $ k \in \{3,4\} $. Уже при $ k = 5 $ мы получаем в точности 
ту же константу $ 1.058... $, полагая $ b = \frac{1}{2} $, $ a = \frac{1}{6} $. Ниже мы приводим таблицу результатов, вытекающих из 
теоремы 4 при различных значениях $ k $. 

\newpage

\begin{longtable}{|c|c|c|c|c|c|c|c|} 
\hline $k$ & $b$ & $a$ & $\zeta_{{\rm clique}}(k) \geq $ & $k$ & $b$ & $a$ & $\zeta_{{\rm clique}}(k) \geq $\\
\hline 
6 &	0.52 &	0.208 &	1.07429 & 20&	0.5&	0.236&	1.12294 \\
\hline
7&	0.498&	0.198&	1.08575 & 25&	0.47&	0.216&	1.12557 \\
\hline
8&	0.5&	0.208&	1.09332 & 30&	0.47&	0.24&	1.12794 \\
\hline
9&	0.5&	0.214&	1.09924 & 35&	0.5&	0.224&	1.12997 \\
\hline
10&	0.5&	0.218&	1.10331 & 40&	0.476&	0.242&	1.13107 \\
\hline
11&	0.5&	0.222&	1.10749 & 45&	0.5&	0.228&	1.13237 \\
\hline
12&	0.5&	0.224&	1.10961 & 50&	0.5&	0.244&	1.13237 \\
\hline
13&	0.474&	0.202&	1.11244 & 55&	0.472&	0.218&	1.13276 \\
\hline
14&	0.472&	0.202&	1.1144 & 101&	0.476&	0.224&	1.13595 \\
\hline
15&	0.5&	0.23&	1.11615 & 1001&	0.484&	0.234&	1.13917 \\
\hline
16&	0.5&	0.232&	1.11839 & 10001& 0.495&	0.245&	1.1397 \\
\hline
\end{longtable}

Отметим еще одно важное обстоятельство. Из теоремы 5 видно, что если $ l = \frac{n}{2} $ и $ \delta = cn $, $ c < \frac{1}{4} $, 
то размер максимального кода (максимальной клики) не превосходит константы, которая тем больше, чем ближе $ c $ к одной четверти. 
Однако хорошо известно, что при $ c = \frac{1}{4} $ клики в соответствующем дистанционном графе имеют растущий размер, а возможно, 
даже размер $ n $. Это связано с классической задачей про матрицы Адамара (см. \cite{Hall}). 

Указанное обстоятельство в некотором смысле послужит основой для применения в следующем разделе вероятностного метода, который для 
большинства значений $ k $ даст лучшие оценки величины $ \zeta_{{\rm clique}}(k) $, нежели те, которые мы уже получили.

\section{Оценка величины $ \zeta_{{\rm clique}}(k) $ с помощью вероятностного метода}

Общая идея подхода состоит в том, чтобы взять графы, у каждого из которых изначально есть клики сколь угодно большого 
размера, а потом удалить часть ребер из этих графов, уничтожив в них все клики данной величины, и выбрать из множества 
рассмотренных графов оптимальный. Наиболее эффективным окажется здесь вероятностный метод (см. \cite{AS}, \cite{Bol}, \cite{Rai7}). 

Заметим, что результаты этого раздела очень похожи на результаты работ \cite{RaiRub1}, \cite{RaiRub2}. Однако будут тут и небольшие, но важные 
отличия, которые приведут к некоторому уточнению оценок. Поэтому и формулировки, и доказательства мы воспроизведем во всех 
подробностях.  

Итак, в параграфе 4.1 мы сформулируем основные результаты. В параграфе 4.2 мы проведем сравнение новых результатов как между собой, так и с 
результатами раздела 3. В параграфах 4.3 и 4.4 мы докажем теоремы из параграфа 4.1. 

\subsection{Формулировки результатов}

Следующие две теоремы дают общий рецепт получения "оптимальных" 
оценок для величин $ \zeta_{{\rm clique}}(k) $. 

\vskip+0.2cm

\noindent {\bf Теорема 6.} {\it Пусть $ k \ge 5 $ -- произвольное натуральное 
число, а $ b_0 \in \left(0,\frac{1}{2}\right) $ -- произвольное 
вещественное число. Положим 
$$
\tau_0 = \tau_0(b_0) = \left(\frac{b_0}{2}\right)^{-\frac{b_0}{2}}
\left(1-\frac{b_0}{2}\right)^{-\left(1-\frac{b_0}{2}\right)}, ~~
\tau_1 = \tau_1(b_0)=b_0^{-b_0} (1-b_0)^{-(1-b_0)}.
$$
Рассмотрим 
$$ 
{\cal C} = {\cal C}(b_0,k)=\left\{c' \in [\tau_0,\tau_1]: ~ 
k \ge \left \lceil \frac{2 \ln \tau_1}{\ln c' - \ln \tau_0} \right \rceil  
\right\}.
$$
Положим $ c = c(b_0,k)=\inf {\cal C} $, если $ {\cal C} \neq \emptyset $, 
и $ c = \tau_1 $ иначе. 
Тогда $ \zeta_{{\rm clique}}(k) \ge \frac{\tau_1}{c} $.}

\vskip+0.2cm

\noindent {\bf Теорема 7.} {\it Пусть $ k \ge 5 $ -- произвольное натуральное 
число, а $ b_{-1}, b_1 $ -- произвольные 
вещественные числа, удовлетворяющие ограничениям
$$
b_{-1}, b_1 \in (0,1), ~~ b_{-1}+b_1 \le \frac{1}{2}, ~~
b_{-1} \le b_1.
$$
Положим
$$
A = \frac{2+9b_{-1}+3b_1-\sqrt{(2+9b_{-1}+3b_1)^2-12(3b_{-1}+
b_1)^2}}{12},
$$
$$
B = \frac{3b_{-1}+b_1}{2}-2A, \, \, \,
C = 1 + A - \frac{3b_{-1}+b_1}{2}.
$$
Пусть, далее, 
$$
\rho_0 = \rho_0(b_{-1},b_1)=A^{-A} B^{-B} C^{-C}, ~~ 
\rho_1 = \rho_1(b_{-1},b_1)
= b_{-1}^{-b_{-1}} b_1^{-b_1} (1-b_{-1}-b_1)^{-(1-b_{-1}-b_1)}.
$$
Рассмотрим 
$$ 
{\cal C} = {\cal C}(b_{-1},b_1,k)=\left\{c' \in [\rho_0,\rho_1]: ~ 
k \ge \left \lceil \frac{2 \ln \rho_1}{\ln c' - \ln \rho_0} \right \rceil  
\right\}.
$$
Положим $ c = c(b_{-1},b_1,k)=\inf {\cal C} $, 
если $ {\cal C} \neq \emptyset $, и $ c = \rho_1 $ иначе. 
Тогда $ \zeta_{{\rm clique}}(k) \ge \frac{\rho_1}{c} $.}
\vskip+0.2cm

Нетрудно проверить, что в условиях теоремы 6 величина $ \tau_0 $ всегда 
меньше величины $ \tau_1 $, а потому множество $ {\cal C} $ 
определено корректно. В конечном счете имеем
$$
\zeta_{{\rm clique}}(k) \ge \max\limits_{b_0 \in \left(0,\frac{1}{2}\right)} 
\frac{\tau_1(b_0)}{c(b_0,k)}.
$$

Аналогичный комментарий верен и относительно теоремы 7. Правда, там не 
столь очевидно неравенство $ \rho_0 < \rho_1 $; доказательство этого 
неравенства можно найти в \cite{Rai4}. В итоге 
$$
\zeta_{{\rm clique}}(k) \ge \max\limits_{b_{-1}, b_1} 
\frac{\rho_1(b_{-1},b_1)}{c(b_{-1},b_1,k)}.
$$

В следующем параграфе мы сравним результаты теорем 4, 6 и 7.

\subsection{Сопоставление результатов}

Прежде всего заметим, что теоремы 6 и 7 отличаются от своих аналогов из \cite{RaiRub2} тем, что 
в \cite{RaiRub2} величина $ \zeta(k) $ определялась условием $ \omega(G) \le k $, тогда как теперь
мы делаем более сильное ограничение $ \omega(G) < k $. Тем не менее, дополнительный вероятностный трюк, 
который мы ввели в доказательства, позволил сохранить результаты в точности такими же, какими они были прежде -- 
при нестрогой оценке размера максимальной клики. Соответственно, мы собрали оптимальные 
результаты в таблице, которую приводим ниже. В первом столбце указано $ k $, 
во втором и третьем -- $ b_{-1}, b_1 $ (параметры из теоремы 7), 
в четвертом -- наилучшая оценка в 
теореме 7, в пятом -- $ b_0 $ (параметр из теоремы 6), в шестом -- 
наилучшая оценка в теореме 6, в седьмом -- лучшая из оценок теорем 6 и 7; при  
этом в четвертом и шестом столбцах жирным шрифтом выделена та оценка, которая 
попала в седьмой столбец. Все оценки даны с восемью точными знаками после 
запятой.

\begin{longtable}{|c|c|c|c|c|c|c|} \hline
 & \multicolumn{3}{|c|}{Теорема 7} & \multicolumn{2}{|c|}{Теорема 6} & \\
\hline $k$ & $b_{-1}$ & $b_1$ & $\zeta(k)\geq$ & $b_0$ &
$\zeta(k)\geq$ & $\zeta(k)\geq$ \endhead
\hline 5 & 0.00001734 & 0.02544379 & 1.00280720 & 0.02825291 & \textbf{1.00297305}&  1.00297305 \\
\hline 6 & 0.00071800 & 0.08438940 & 1.01694361 & 0.09659179 & \textbf{1.01838640}&  1.01838640 \\
\hline 7 & 0.00270696 & 0.12703569 & 1.03395854 & 0.14157882 & \textbf{1.03626539}&  1.03626539 \\
\hline 8 & 0.00553133 & 0.15711445 & 1.05011426 & 0.17051644 & \textbf{1.05242870}&  1.05242870 \\
\hline 9 & 0.00869977 & 0.17904575 & 1.06460594 & 0.19035277 & \textbf{1.06632451}&  1.06632451 \\
\hline 10 & 0.01190629 & 0.19555639 & 1.07739346 & 0.20472195 & \textbf{1.07816823}&  1.07816823 \\
\hline 11 & 0.01498920 & 0.20833110 & \textbf{1.08864080} & 0.21558694 & 1.08829570 &  1.08864080 \\
\hline 12 & 0.01787410 & 0.21844807 & \textbf{1.09855184} & 0.22408139 & 1.09701591 &  1.09855184 \\
\hline 13 & 0.02053497 & 0.22662152 & \textbf{1.10731983} & 0.23090096 & 1.10458387 &  1.10731983 \\
\hline 14 & 0.02297102 & 0.23333932 & \textbf{1.11511348} & 0.23649470 & 1.11120346 &  1.11511348 \\
\hline 15 & 0.02519380 & 0.23894375 & \textbf{1.12207552} & 0.24116490 & 1.11703654 &  1.12207552 \\
\hline 16 & 0.02722013 & 0.24368069 & \textbf{1.12832523} & 0.24512232 & 1.12221184 &  1.12832523 \\
\hline 17 & 0.02906839 & 0.24773059 & \textbf{1.13396195} & 0.24851823 & 1.12683247 &  1.13396195 \\
\hline 18 & 0.03075664 & 0.25122837 & \textbf{1.13906849} & 0.25146404 & 1.13098160 &  1.13906849 \\
\hline 19 & 0.03230179 & 0.25427658 & \textbf{1.14371408} & 0.25404355 & 1.13472693 &  1.14371408 \\
\hline 20 & 0.03371912 & 0.25695447 & \textbf{1.14795689} & 0.25632099 & 1.13812400 &  1.14795689 \\
\hline 25 & 0.03931602 & 0.26656957 & \textbf{1.16461196} & 0.26461281 & 1.15125524 &  1.16461196 \\
\hline 30 & 0.04321023 & 0.27249456 & \textbf{1.17618197} & 0.26984084 & 1.16020045 &  1.17618197 \\
\hline 35 & 0.04605883 & 0.27649512 & \textbf{1.18467380} & 0.27343793 & 1.16668123 &  1.18467380 \\
\hline 40 & 0.04822699 & 0.27937171 & \textbf{1.19116666} & 0.27606415 & 1.17159101 &  1.19116666 \\
\hline 45 & 0.04992995 & 0.28153709 & \textbf{1.19628994} & 0.27806578 & 1.17543859 &  1.19628994 \\
\hline 50 & 0.05130175 & 0.28322476 & \textbf{1.20043462} & 0.27964192 & 1.17853468 &  1.20043462 \\
\hline 55 & 0.05242980 & 0.28457649 & \textbf{1.20385612} & 0.28091520 & 1.18107971 &  1.20385612 \\
\hline 100 & 0.05755871 & 0.29034354 & \textbf{1.21959032} & 0.28647009 & 1.19266305 &  1.21959032 \\
\hline 1000 & 0.06326998 & 0.29611686 & \textbf{1.23753236} & 0.29226811 & 1.20564840 &  1.23753236 \\
\hline 10000 & 0.06384354 & 0.29666283 & \textbf{1.23936294} & 0.29283085 & 1.20696080 &  1.23936294 \\
\hline 100000 & 0.06390091 & 0.29671712 & \textbf{1.23954636} & 0.29288697 & 1.20709218 &  1.23954636 \\
\hline \end{longtable}

Видно, что теорема 6 сильнее теоремы 7 лишь при $ k \le 10 $. Любопытно, что теорема 6 доказывается с помощью 
(0,1)-векторов, а теорема 7 -- с помощью (-1,0,1)-векторов. Иными словами, если в подходе с равновесными кодами (вернее, 
в его "ручном" аналоге из раздела 2) добавление минус единиц не помогало (см. \S 2.4), то в вероятностном подходе 
увеличение числа координат дает уточнения результатов. 

Далее, число 
в шестой колонке, очевидно, стремится к 1.207..., а число в четвертой колонке становится все ближе к 1.239... Это 
и суть константы из работ \cite{FW} и \cite{Rai2}, в которых получены самые сильные из известных нижние оценки для 
$ \chi({\mathbb R}^n) $ (см. раздел 1). Таким образом, в асимптотике по $ k $ теоремы 6 и 7 дают очень сильные результаты. 

Результат теоремы 4 в указанном смысле слабее. Оценки в таблице из \S 3.3 не стремятся ни к 1.207, ни к 1.239. Их 
пределом является величина 
$$
1.139... = 2 \cdot \left(\frac{1}{4}\right)^{1/4} \cdot \left(\frac{3}{4}\right)^{3/4}. 
$$
Это видно из таблицы, но это нетрудно доказать формально. 

В то же время при малых $ k $ теорема 4 все-таки лучше. А именно, она сильнее теорем 6 и 7 при всех $ k \le 13 $. 
Таким образом, результат теоремы 6 оказывается полностью в некотором смысле покрыт результатом теоремы 4. 

Заметим, наконец, что теоремы 6 и 7, по существу, не применимы к $ k \le 4 $. Таким образом, в этой области "штучные подходы" и их 
обобщения через "кодирование" просто незаменимы. 

В следующем разделе мы докажем теорему 6. 

\subsection{Доказательство теоремы 6}

Зафиксируем $ k, b_0 $, а стало быть, и $ \tau_0, \tau_1, {\cal C}, c $.
Легко видеть, что либо $ {\cal C} = \emptyset $, либо $ {\cal C} = \{\tau_1\} $,
либо $ {\cal C} = [c,\tau_1] $, $ c \ge \tau_0 $. В первых двух случаях 
утверждение теоремы тривиально. Рассмотрим ситуацию, когда 
$ \tau_0 \le c < \tau_1 $. 

Возьмем произвольное $ c' \in (c,\tau_1) $. Здесь важно, что $ c' $ строго 
больше $ c $, хотя и сколь угодно близко к $ c $. Если мы покажем, 
что $ \zeta_{{\rm clique}}(k) \ge \frac{\tau_1}{c'} $, то, беря инфимум обеих частей 
неравенства по $ c' $, мы получим искомое утверждение    
$ \zeta_{{\rm clique}}(k) \ge \frac{\tau_1}{c} $. 

Итак, нам нужно убедиться в существовании такой
функции $ \delta(n) = o(1) $, что при всех $ n $ найдется граф расстояний 
$ G=(V,E) $ в $ {\mathbb R}^n $, у которого одновременно $ \omega(G) < k $
и $ \chi(G) \ge \left(\frac{\tau_1}{c'}+\delta(n)\right)^n $. 

Положим для каждого достаточно большого $ n \in {\mathbb N} $
$$
{\cal V}_n = \{{\bf x} = (x_1, \dots, x_n): ~ x_i \in \{0,1\}, ~ 
|\{i:~x_i=1\}| = [b_0n]\}.
$$
Пусть $ p $ -- минимальное простое число, удовлетворяющее условию 
$ [b_0n] - 2p < 0 $. Положим 
$$
{\cal E}_n = \left\{\{{\bf x},{\bf y}\}:~ |{\bf x}-{\bf y}|=\sqrt{2p}\right\}, ~~ 
{\cal G}_n = ({\cal V}_n,{\cal E}_n).
$$

С помощью формулы Стирлинга легко показать, что 
$$ 
N = |{\cal V}_n| = 
(\tau_1+o(1))^n. 
$$ 
В то же время уже хорошо известный нам линейно-алгебраический метод 
дает оценку 
$$ 
\alpha = 
\alpha({\cal G}_n) \le n C_n^{p}. 
$$ 
При $ n \to \infty $ имеем 
$ p \sim \frac{b_0}{2} n $, и, значит, за счет 
все той же формулы Стирлинга $ \alpha \le (\tau_0+\delta_1)^n $
с некоторой $ \delta_1 = o(1) $.                                   

Конечно, если снова вспомнить простую оценку $ \chi(G) \ge \frac{|V|}{\alpha(G)} $,
то окажется, что 
$$ 
\chi({\cal G}_n) \ge \left(\frac{\tau_1}{\tau_0}+\delta_2(n)\right)^n 
$$
с некоторой $ \delta_2 = o(1) $, и это, на 
первый взгляд, даже лучше заявленного в теореме результата. Неприятность в 
том, что граф расстояний $ {\cal G}_n $ вполне может содержать 
(и содержит) клики 
размера больше $ k $ (ср. конец \S 3.3). 
 
Воспользуемся вероятностным методом (см. \cite{AS}, \cite{Bol}, \cite{Rai7}). 
Рассмотрим произвольное число $ \gamma \in \left(\frac{\tau_0}{c'},1\right) $.
Такое $ \gamma $ существует, поскольку $ c' > c \ge \tau_0 $. Положим 
$ q = \gamma^n $, где $ n $ по-прежнему достаточно велико. 
Построим случайный подграф $ G = ({\cal V}_n,E) $ графа $ {\cal G}_n $, 
"кладя" в $ E $ каждое ребро из $ {\cal E}_n $ с вероятностью $ q $ 
независимо от остальных ребер. Образуется вероятностное 
пространство $ (\Omega_n,{\cal B}_n,P_n) $, в котором 
$$ 
\Omega_n = 
\{G=({\cal V}_n,E), ~ E \subseteq {\cal E}_n\}, ~~ 
{\cal B}_n = 2^{\Omega_n}, ~~ P_n(G) = q^{|E|}(1-q)^{|{\cal E}_n|-|E|} ~~~
{\rm для} ~~ G = ({\cal V}_n,E).
$$ 

Определим на $ \Omega_n $ два семейства случайных величин. Во-первых,
при каждом $ l \in {\mathbb N} $ и при каждом $ G \in \Omega_n $
положим $ X_l(G) $ равной числу независимых множеств вершин размера $ l $ 
в графе $ G $. Во-вторых, при аналогичных условиях $ Y_m(G) $ -- это 
число клик размера $ m $ в $ G $. 

Пусть $ l = \left[(c')^n\right] $. Ясно, что, за счет неравенства 
$ c' < \tau_1 $, при больших $ n $ выполнено $ l < N = |{\cal V}_n| $, а стало 
быть, величина $ X_l $ корректно определена. 

Допустим, мы показали, что  
$$
MX_l < \frac{1}{2}, ~~ MY_{k} < \frac{N}{4},
$$
где через $ M\xi $ обозначено математическое ожидание случайной величины $ \xi $. 
Тогда по неравенству Маркова (Чебышёва) имеют место оценки (см. \cite{Gned})
$$
P_n(X_l \ge 1) \le MX_l < \frac{1}{2}, ~~~
P_n\left(Y_k \ge \frac{N}{2}\right) \le \frac{MY_k}{N/2} < \frac{1}{2}, 
$$
т.е. 
$$
P_n(X_l = 0) > \frac{1}{2}, ~~~
P_n\left(Y_k < \frac{N}{2}\right) > \frac{1}{2}. 
$$
Значит, существует граф $ G $ в $ {\mathbb R}^n $ с 
$ \alpha(G) \le l $ и с не более, чем $ \frac{N}{2} $, кликами размера $ k $. Удалим 
из каждой клики в $ G $ по одной вершине, получится граф $ G' $ с не менее $ \frac{N}{2} $ 
вершинами и с $ \alpha(G') \le l $. В итоге
$$ 
\chi(G') \ge 
\frac{N/2}{\alpha(G')} = \left(\frac{\tau_1}{c'} + \delta(n)\right)^n,
$$
и теорема доказана. (Cледует только заметить, что при малых $ n $ 
доказывать нечего, ведь мы вольны на начальном отрезке натурального 
ряда выбрать значения $ \delta(n) $ сколь угодно большими по модулю.)
Что ж, будем оценивать математические ожидания. 

Начнем с $ MX_l $. За счет линейности математического ожидания, имеем 
$$
MX_l = \sum_{W \subset {\cal V}_n, ~ |W| = l} (1-q)^{r(W)}, ~ 
r(W) = |\{\{{\bf x},{\bf y}\} \in {\cal E}_n: ~ {\bf x} \in W, ~
{\bf y} \in W\}|.
$$

Ясно, что, поскольку $ c' > \tau_0 $, то $ l > \alpha = \alpha({\cal G}_n) $
при всех достаточно больших $ n $, и, стало быть, для каждого $ W \subset 
{\cal V}_n $, $ |W| = l $, выполнено $ r(W) > 0 $. Нетрудно, однако, 
установить и гораздо более точное неравенство $ r(W) \ge \frac{l^2}{4 \alpha} $, верное
при $ l > \alpha $ (см. \cite{RaiRub2}).

Заметим, что 
$$
\frac{l^2}{4\alpha} \ge \frac{((c')^2+\delta_3(n))^n}{(\tau_0+\delta_4(n))^n} = 
\left(\frac{(c')^2}{\tau_0}+\delta_5(n)\right)^n, ~~ \delta_3(n) = o(1), ~~
\delta_4(n)=o(1), ~~ \delta_5(n) = o(1).
$$
Таким образом, полагая 
$$
A_l = \frac{l^2}{4\alpha}, 
$$
имеем (с некоторыми $ \delta_i(n) = o(1) $)
$$
MX_l \le C_N^l \cdot (1-q)^{A_l} \le \left(\frac{eN}{l}\right)^l \cdot
(1-q)^{A_l} \le 
$$
$$
\le \left(\frac{\tau_1}{c'}+\delta_6(n)\right)^{(c'+\delta_7(n))^n}
\cdot e^{-q\left(\frac{(c')^2}{\tau_0}+\delta_5(n)\right)^n} = 
e^{(c'+\delta_8(n))^n - \left(\frac{\gamma \cdot (c')^2}{\tau_0}+
\delta_9(n)\right)^n}.
$$

У нас $ \gamma > \frac{\tau_0}{c'} $. Значит, $ \frac{\gamma \cdot 
(c')^2}{\tau_0} > c' $, т.е.     
$$
(c'+\delta_8(n))^n - \left(\frac{\gamma \cdot (c')^2}{\tau_0}+
\delta_9(n)\right)^n \to -\infty,
$$
а 
$$
e^{(c'+\delta_8(n))^n - \left(\frac{\gamma \cdot (c')^2}{\tau_0}+
\delta_9(n)\right)^n} \to 0.
$$
Следовательно, при всех достаточно больших $ n $ выполнено $ MX_l < 
\frac{1}{2} $, и нам остается обосновать оценку $ MY_{k} < \frac{N}{2} $.

За счет линейности математического ожидания имеем
$ MY_{k} \le C_N^{k} q^{C_{k}^2} $. Далее, поскольку $ k $ -- константа
(т.е. $ k $ не зависит от $ n $), 
$$
MY_{k} \le N^{k} \cdot \left(\gamma^n\right)^{C_{k}^2} =  
(\tau_1+\delta_{10}(n))^{nk} \gamma^{n \frac{k(k-1)}{2}} = 
(\tau_1+\delta_{10}(n))^n e^{n(k-1)(\ln \tau_1 + \frac{k}{2} \ln \gamma + \delta_{11}(n))}.
$$

Из определения $ c $ и из условия $ c < c' $ видно, что 
$$
k \ge \left \lceil \frac{2 \ln \tau_1}{\ln c - \ln \tau_0} \right \rceil 
\ge \frac{2 \ln \tau_1}{\ln c - \ln \tau_0} >  
\frac{2 \ln \tau_1}{\ln c' - \ln \tau_0}. 
$$

Все наши прежние выкладки верны при любом $ \gamma > \frac{\tau_0}{c'} $.
Возьмем $ \gamma $ настолько близким к своей нижней границе, чтобы 
оказалась справедливой оценка $ k > \frac{2 \ln \tau_1}{-\ln \gamma} $. 
Тогда, c учетом отрицательности величины $ \ln \gamma $, 
$$
\ln \tau_1 + \frac{k}{2} \ln \gamma < \ln \tau_1 - \ln \tau_1 = 0.
$$
Таким образом, $ MY_{k} = o(N) $ при $ n \to \infty $, и нужное нам 
неравенство при больших $ n $ имеет место.

Теорема доказана.

\subsection{Доказательство теоремы 7}

Схема действий практически полностью повторяет свой аналог из предыдущего 
параграфа. По существу, разница только в построении графа $ {\cal G}_n $. 

Пусть
$$
{\cal V}_n = \{{\bf x} = (x_1, \dots, x_n): ~ x_i \in \{-1,0,1\}, ~ 
|\{i:~x_i=-1\}| = [b_{-1}n], ~ |\{i:~x_i=1\}| = [b_{1}n]\}.
$$
Пусть, далее, $ p $ -- минимальное простое число, удовлетворяющее условию 
$ [b_{-1}n] + [b_1n] - 2p < -2[b_{-1}n] $. Положим 
$$
{\cal E}_n = \left\{({\bf x},{\bf y}):~ |{\bf x}-{\bf y}|=\sqrt{2p}\right\}, ~~ 
{\cal G}_n = ({\cal V}_n,{\cal E}_n).
$$

По Стирлингу, $ |{\cal V}_n| = (\rho_1+o(1))^n $. Кроме того, известно
(см. \cite{Rai4}), что 
$$ 
\alpha=\alpha({\cal G}_n) \le (\rho_0+o(1))^n.
$$ 
Наконец, неравенство типа $ r(W) \ge \frac{l^2}{4\alpha} $ верно и здесь.
Дальнейшие рассуждения ясны, и теорема доказана.

\section{О величине $ \zeta_{{\rm girth}}(k) $}

Этот раздел мы разобьем на три части. В первой из них мы сформулируем наши основные результаты. 
Во второй и третьей частях мы эти результаты докажем. 

\subsection{Формулировки результатов}

Для получения оценок величины $ \zeta_{{\rm clique}}(k) $ мы использовали графы $ G_{b,a} = (V_b,E_a) $ следующего вида: 
$$
V_b = \{{\bf x} = (x_1, \dots, x_n): ~ \forall ~i~~ x_i \in \{0,1\}, ~ x_1 + \ldots + x_n = b\}, ~~
E_a = \{\{{\bf x},{\bf y}\}:~ {\bf x}, {\bf y} \in V, ~ ({\bf x},{\bf y}) = a\}
$$
при некоторых $ b $ и $ a $. Оказывается, в таких графах всегда есть четные циклы любой допустимой длины. В частности, 
справедлива 

\vskip+0.2cm

\noindent {\bf Теорема 8.} {\it Пусть $ b \le \frac{n}{2} $, $ a \ge 1 $, $ k \in {\mathbb N} $, причем $ k \ge 2 $ и $ b - k \ge a $. Тогда в графе 
$ G_{b,a} $ есть цикл длины $ 2k $.}

\vskip+0.2cm

Теорему 8 мы докажем в следующем параграфе. Уже сейчас, однако, ясно, что величину $ \zeta_{{\rm girth}}(k) $ не удастся оценить 
с помощью какого бы то ни было графа $ G_{b,a} $. Графы с более сложной координатной структурой вершин (например, графы, у которых 
вершины могут содержать -1, а не только 0 и 1) довольно трудно исследовать на предмет наличия в них каких-либо циклов, и мы их 
также в настоящей работе не рассматриваем. Зато в рамках нашей стандартной технологии мы найдем графы вида $ G_{b,a} $, которые 
не содержат {\it нечетных} циклов длины $ \le 2k+1 $ с любым наперед заданным $ k $ и которые, тем не менее, имеют экспоненциально 
большие хроматические числа. Это весьма любопытно, ведь, как известно, критерием двудольности графа служит отсутствие нечетных циклов в нем. 
Конечно, мы избавимся лишь от коротких циклов, и все же это по-прежнему нетривиально. 

Итак, обозначим через $ g_{{\rm odd}}(G) $ длину кратчайшего нечетного цикла в $ G $. Положим
$$
\zeta_{{\rm odd ~ girth}}(k) = \sup \left\{\zeta: ~ \exists ~ {\rm функция} ~ \delta=\delta(n), ~ {\rm такая}, ~ {\rm что} ~ 
\lim\limits_{n \to \infty} \delta(n) = 0 \right.
$$
$$
\left. \phantom{\lim\limits_{n \to \infty}} {\rm и} ~ \forall ~ n ~
\exists ~ G ~ {\rm в} ~ {\mathbb R}^n, ~ {\rm у} ~ {\rm которого} ~ g_{{\rm odd}}(G) > k, ~ 
\chi(G) \ge (\zeta+\delta(n))^n\right\}.
$$

Справедлива 
\vskip+0.2cm

\noindent {\bf Теорема 9.} {\it Имеет место оценка
$$
\zeta_{{\rm odd ~ girth}}(2k+1) \ge 2 \cdot \left(\frac{k}{2k+1}\right)^{\frac{k}{2k+1}} \cdot \left(\frac{k+1}{2k+1}\right)^{\frac{k+1}{2k+1}}.
$$}

\vskip+0.2cm

Теорему 9 мы докажем в параграфе 5.3. Стоит отметить, что 
$$ 
\zeta_{{\rm odd ~ girth}}(3) = \zeta_{{\rm clique}}(3).
$$
Нетрудно видеть, что в указанном случае оценки из теорем 3 и 9 совпадают.

\subsection{Доказательство теоремы 8}

Нам нужно найти в графе $ G_{b,a} $ цикл длины $ 2k $. Иными словами, нам надо предъявить вершины $ {\bf x}_1, \dots, {\bf x}_{2k} \in V_b $, для 
которых $ ({\bf x}_i,{\bf x}_{i+1}) = a $ при всех $ i \in \{1, \dots, 2k\} $ (здесь $ 2k+1 = 1 $). Сейчас мы эти вершины просто укажем. Для 
этого сперва разобьем множество $ \{1, \dots, n\} $ координатных позиций каждого вектора из $ V_b $ на части $ A_1, \dots, A_5 $:
$$
A_1 = \{1, \dots, a-1\}, ~ A_2 = \{a, \dots, b-k\}, ~ A_3 = \{b-k+1, \dots, 2b-2k-a+1\}, ~ 
$$
$$
A_4 = \{2b-2k-a+2, \dots, 2b-a+1\}.
$$
Ясно, что 
$$
|A_1| = a-1, ~ |A_2| = b-a-k+1, ~ |A_3| = b-a-k+1, ~ |A_4| = 2k.
$$
При этом все корректно, поскольку 
$$
b-k \ge a, ~~~ 2b - a + 1 \le 2 \cdot \frac{n}{2} - 1 + 1 = n,
$$
так что $ |A_5| \ge 0 $. 

У каждого из будущих векторов $ {\bf x}_1, \dots, {\bf x}_{2k} $ в части $ A_1 $ мы расположим единицы, в части $ A_5 $ -- нули. Далее, наполнение 
частей $ A_2, A_3 $, имеющих одинаковую мощность, мы будем чередовать. В векторах с нечетными номерами мы $ A_2 $ заполним единицами, $ A_3 $ -- 
нулями, а в векторах с четными номерами мы, напротив, $ A_2 $ заполним нулями, $ A_3 $ -- единицами. Наконец, с частью $ A_4 $ мы поступим 
следующим образом. В векторе $ {\bf x}_1 $ на первые $ k $ позиций из $ A_4 $ поставим единицы, на остальные позиции поставим нули. В векторе 
$ {\bf x}_{i+1} $ мы поставим 1 на последнюю позицию, в которой у вектора $ {\bf x}_i $ стоит 1, а за ней последовательно разместим еще $ k-1 $ 
единиц (мы считаем, что $ 2k+t = t $). Теперь легко видеть, что 
$$ 
({\bf x}_{i+1}, {\bf x}_i) = a-1 + 1 = a,
$$
так что вершины с соседними номерами соединены ребрами, причем часть $ A_4 $ устроена так, что указанное соотношение верно и при $ i = 2k $. 
Кроме того, число единиц в каждом векторе равно 
$$
(a-1) + (b-a-k+1) + 2k = b.
$$
Значит, векторы $ {\bf x}_1, \dots, {\bf x}_{2k} $ образуют цикл длины $ 2k $ в графе $ G_{b,a} $, и теорема доказана. 

\subsection{Доказательство теоремы 9}

Положим $ b = \left[\frac{n}{2}\right] $. Возьмем минимальное простое число $ p $, большее величины $ \frac{nk}{2k+1} $. Ясно, что, как всегда, 
$ p \sim \frac{nk}{2k+1} $ и что при достаточно больших $ n $ выполнены неравенства 
$$ 
b > p > \frac{nk}{2k+1}. 
$$
Пусть $ a = b-p $. Очевидно, $ b - 2p < 0 $. Сказанное сразу же влечет оценку 
$$
\chi(G_{b,a}) \ge \frac{C_n^b}{n C_n^p}, 
$$
которая получается за счет линейной алгебры. В свою очередь, стандартная аналитика (формула Стирлинга) показывает, что 
$$
\frac{C_n^b}{n C_n^p} = \left(2 \cdot \left(\frac{k}{2k+1}\right)^{\frac{k}{2k+1}} \cdot \left(\frac{k+1}{2k+1}\right)^{\frac{k+1}{2k+1}}+o(1)\right)^n.
$$
Таким образом, остается убедиться в том, что $ g_{{\rm odd}}(G_{b,a}) > 2k+1 $.  

Предположим, в $ G_{b,a} $ есть цикл на вершинах $ {\bf x}_1, \dots, {\bf x}_{2k+1} $. Обозначим через $ A_i $, $ i = 1, \dots, 2k+1 $, 
множество координатных позиций, на которых в векторе $ {\bf x}_i $ стоят единицы. 

\vskip+0.2cm

\noindent{\bf Лемма 1.} {\it Пусть $ |A_1 \cap A_i| = \alpha $. Тогда $ |A_1 \cap A_{i+2}| \le \alpha + n + 2a - 2b $.} 

\vskip+0.2cm

\paragraph{Доказательство леммы 1.} Посмотрим сперва на вектор $ {\bf x}_{i+1} $ и соответствующее множество $ A_{i+1} $. 
Поскольку 
$$
|\{1, \dots, n\} \setminus (A_1 \cup A_i)| = n-2b+\alpha, 
$$
количество элементов множества $ A_{i+1} $, находящихся в множестве $ \{1, \dots, n\} \setminus (A_1 \cup A_i) $, 
можно представить в виде $ n-2b+\alpha-\varepsilon $, $ \varepsilon \ge 0 $. При этом, конечно, $ |A_i \cap A_{i+1}| = a $. Значит, 
$$
|A_{i+1} \cap (A_1 \setminus A_i)| = b - a - (n-2b+\alpha-\varepsilon) = 3b-a-\alpha-n+\varepsilon. 
$$
Положим 
$$ 
|A_{i+1} \cap (A_i \cap A_1)| = a - \delta, ~~ |A_{i+1} \cap (A_i \setminus A_1)| = \delta, ~~ 0 \le \delta \le a.                    
$$

Теперь рассмотрим $ A_{i+2} $. Дабы понять, в каком случае пересечение этого множества с $ A_1 $ имеет наибольшую мощность, разобьем 
$ A_1 $ на части: 
$$
A_1 = (A_{i+1} \cap (A_1 \setminus A_i)) \bigsqcup (A_{i+1} \cap (A_i \cap A_1)) \bigsqcup B_3 = B_1 \bigsqcup B_2 \bigsqcup B_3.
$$
Величина $ |A_{i+2} \cap A_1| $ максимальна, коль скоро $ A_{i+2} \cap A_{i+1} \subseteq B_1 \cup B_2 $. Следовательно, 
$$
|A_{i+2} \cap A_1| \le |A_{i+2} \cap A_{i+1}| + |B_3| = a + b - (3b-a-\alpha-n+\varepsilon) - (a-\delta) = n - 2b + a + \alpha - 
\varepsilon + \delta \le n - 2b + 2a + \alpha.
$$

Лемма доказана. 
 
\vskip+0.2cm

Мы знаем, что $ |A_1 \cap A_2| = a $. Стало быть, лемма 1 говорит, что 
$$ 
|A_1 \cap A_4| \le a + (n+2a-2b), ~~ |A_1 \cap A_6| \le a + 2(n+2a-2b), ~~, 
|A_1 \cap A_{2k+2}| \le a + k(n+2a-2b), 
$$
где $ A_{2k+2} = A_1 $. Если мы сейчас докажем, что $ a+k(n+2a-2b) < b $, то мы, очевидно, придем к противоречию, ведь 
$$ 
|A_1 \cap A_{2k+2}| = |A_1 \cap A_1| = b. 
$$
Итак, 
$$
a+k(n+2a-2b) = b-p + k (n+2(b-p)-2b) = b - p(2k+1) + kn < b - \frac{kn}{2k+1} \cdot (2k+1) + kn = b, 
$$
и теорема 9 доказана.

\section{Несколько запрещенных расстояний}

В этом разделе будет шесть параграфов. В первом из них мы расскажем об истории задачи о хроматическом числе 
пространства с несколькими запрещенными расстояниями и поставим нашу основную задачу. Во втором параграфе мы приведем формулировку основного 
результата. В третьем параграфе мы этот результат докажем. В четвертом параграфе мы дадим несколько комментариев о связи доказанного нами 
результата с его аналогами из предыдущих разделов. Поскольку в формулировке основного результата будет присутствовать оптимизация по 
нескольким параметрам, мы в пятом параграфе расскажем о решении соответствующей оптимизационной задачи (которое само по себе нетривиально). 
Наконец, в шестом параграфе мы приведем некоторые таблицы с численными результатами, которые получаются при решении 
упомянутой выше оптимизационной задачи.  

\subsection{История вопроса и постановка основной задачи}

В настоящем разделе нас будут интересовать дистанционные графы, у которых вершины по-прежнему являются элементами пространства 
$ {\mathbb R}^n $ с евклидовой метрикой, а ребра задаются не одним запрещенными расстоянием, но некоторым конечным множеством 
запретов $ {\cal A} $. А именно, мы прежде всего рассмотрим величину 
$$
\chi({\mathbb R}^n, {\cal A}) = 
\min \left\{\chi:~ {\mathbb R}^n = V_1 \bigsqcup \ldots \bigsqcup V_{\chi}, ~ \forall ~i ~\forall ~ {\bf x}, {\bf y} \in V_i ~~
|{\bf x}-{\bf y}| \not \in {\cal A}\right\}.
$$

Конечно, можно пытаться исследовать различные свойства величины $ \chi({\mathbb R}^n, {\cal A}) $ в зависимости от тех или иных 
характеристик множества $ {\cal A} $. Однако здесь мы ограничимся вариантом задачи, при котором искомой является величина 
$$
\overline{\chi}({\mathbb R}^n;m) = \max_{{\cal A}, \, |{\cal A}|=m} \chi({\mathbb R}^n, {\cal A}).
$$

Величину $ \overline{\chi}({\mathbb R}^n;m) $ рассматривал еще П. Эрдеш, который знал, например, что 
$$
c_1 m \sqrt{\ln m} \le \overline{\chi}({\mathbb R}^2;m) \le c_2 m^2, ~~~ c_1, c_2 > 0.
$$
Подобные результаты описаны в обзоре \cite{Szek}. 

Аккуратно подсчитанные нижние оценки для хроматических чисел $ \overline{\chi}({\mathbb R}^n;m) $, коль скоро $ n \to \infty $, 
впервые были приведены в статье \cite{Rai3}. Там было доказано, что 
$$
\overline{\chi}({\mathbb R}^n;1) \ge (1.239...+o(1))^n, ~~~ \overline{\chi}({\mathbb R}^n;2) \ge (1.439...+o(1))^n.
$$
Кроме того, там была установлена общая нижняя оценка, справедливая даже в случае любой функции $ m = m(n) $: 
$$
\overline{\chi}({\mathbb R}^n;m) \ge (c_1m)^{c_2n}, ~~~~ c_1, c_2 > 0. 
$$
Константы $ c_1, c_2 $ были черезвычайно маленькими (порядка $ 2^{-10000} $), поэтому позже указанный результат уточнялся 
(см. \cite{Rai8}). А в работах \cite{Sh1}, \cite{Sh2} была развита оптимизационная техника, позволившая найти оптимальные 
(в рамках некоторого метода) константы типа 1.239 и 1.439. Перечислим эти результаты: 
$$
\overline{\chi}({\mathbb R}^n;2) \ge (1.465...+o(1))^n, ~~~ \overline{\chi}({\mathbb R}^n;3) \ge (1.667...+o(1))^n, ~~~ 
\overline{\chi}({\mathbb R}^n;4) \ge (1.848...+o(1))^n, 
$$
$$
\overline{\chi}({\mathbb R}^n;5) \ge (2.013...+o(1))^n, ~~~ \overline{\chi}({\mathbb R}^n;6) \ge (2.165...+o(1))^n, ~~~ 
\overline{\chi}({\mathbb R}^n;7) \ge (2.308...+o(1))^n, 
$$
$$
\overline{\chi}({\mathbb R}^n;8) \ge (2.442...+o(1))^n, ~~~ \overline{\chi}({\mathbb R}^n;9) \ge (2.570...+o(1))^n, ~~~ 
\overline{\chi}({\mathbb R}^n;10) \ge (2.691...+o(1))^n, 
$$
$$
\overline{\chi}({\mathbb R}^n;11) \ge (2.807...+o(1))^n, ~~~ \overline{\chi}({\mathbb R}^n;12) \ge (2.919...+o(1))^n, ~~~ 
\overline{\chi}({\mathbb R}^n;13) \ge (3.026...+o(1))^n, 
$$
$$
\overline{\chi}({\mathbb R}^n;14) \ge (3.130...+o(1))^n, ~~~ \overline{\chi}({\mathbb R}^n;15) \ge (3.231...+o(1))^n, ~~~ 
\overline{\chi}({\mathbb R}^n;16) \ge (3.328...+o(1))^n.
$$
Список можно продолжать и дальше. 

Что касается верхних оценок величин $ \overline{\chi}({\mathbb R}^n;m) $, то ничего лучшего не придумано, нежели оценка 
$$
\overline{\chi}({\mathbb R}^n;m) \le (3+o(1))^{nm},
$$
мгновенно вытекающая из неравенства $ \overline{\chi}({\mathbb R}^n;1) \le (3+o(1))^n $, которое мы уже цитировали в разделе 1. 
Понятно, что с ростом $ m $ эта оценка катастрофически удаляется от имеющихся нижних оценок. Впрочем, по-видимому, именно 
верхние оценки подлежат дальнейшему усилению. 

Как и в предшествующей части работы, нас будут интересовать графы расстояний с ограниченным кликовым числом. При этом теперь 
наши графы расстояний -- это подграфы в графах $ \mathfrak{G}_{{\cal A}} = (\mathfrak{V}, \mathfrak{E}_{{\cal A}}) $, у которых 
$$
\mathfrak{V} = {\mathbb R}^n, ~~~ \mathfrak{E}_{{\cal A}} = \{\{{\bf x},{\bf y}\}:~ |{\bf x}-{\bf y}| \in {\cal A}\}.
$$

Основным объектом нашего исследования будет величина 
$$
\zeta_{{\rm clique}}(k,m) = \sup \left\{\zeta: ~ \exists ~ {\rm функция} ~ \delta = \delta(n), 
~ {\rm такая}, ~ {\rm что} ~ \lim\limits_{n \to \infty} \delta(n) = 0 \right.
$$
$$
\left. \phantom{\lim\limits_{n \to \infty}} {\rm и} ~ \forall ~ n ~
\exists ~ {\cal A} ~ {\rm мощности} ~ m ~ {\rm и} ~ \exists ~ G ~ {\rm в} ~ \mathfrak{G}_{{\cal A}}, ~ {\rm у} ~ {\rm которого} ~
\omega(G) < k, ~ 
\chi(G) \ge (\zeta+\delta(n))^n\right\}.
$$
Разумеется, $ \zeta_{{\rm clique}}(k,1) = \zeta_{{\rm clique}}(k) $, так что тут мы сконцентрируемся на случаях $ m \ge 2 $. 

В следующем параграфе мы сформулируем наш основной результат. Заметим, что в этом разделе мы не станем говорить про обхват, т.е. 
про величину 
$$
\zeta_{{\rm girth}}(k,m) = \sup \left\{\zeta: ~ \exists ~ {\rm функция} ~ \delta = \delta(n), 
~ {\rm такая}, ~ {\rm что} ~ \lim\limits_{n \to \infty} \delta(n) = 0 \right.
$$
$$
\left. \phantom{\lim\limits_{n \to \infty}} {\rm и} ~ \forall ~ n ~
\exists ~ {\cal A} ~ {\rm мощности} ~ m ~ {\rm и} ~ \exists ~ G ~ {\rm в} ~ \mathfrak{G}_{{\cal A}}, ~ {\rm у} ~ {\rm которого} ~
g(G) > k, ~ 
\chi(G) \ge (\zeta+\delta(n))^n\right\},
$$
которая также, конечно, имеет смысл. 

\subsection{Формулировка результата}

Нам удалось доказать следующую теорему. 

\vskip+0.2cm

\noindent {\bf Теорема 10.} {\it Пусть даны числа $ k $ и $ m $. Пусть также $ r \ge 1 $ -- произвольное натуральное число. 
Рассмотрим симплекс 
$$
\Delta = \{{\bf v} = (v_0, v_1, \dots, v_r): ~ v_i \in [0,1], ~ v_0 + v_1 + \ldots + v_r = 1\}.
$$
Для каждого $ {\bf v} \in \Delta $ положим 
$$
g(t)=g(t,{\bf v})=\sum_{i=0}^r \mathbb{I}_{\left\{t\geq \sum_{j=0}^i v_j\right\}}(t), ~~~ \overline{s}'=\overline{s}'({\bf v}) = \int_0^1 g^2(t) dt, ~~~
\underline{s}' = \underline{s}'({\bf v}) = \int_0^1 g(t)g(1-t) dt, 
$$
$$
{\bf b} = (0,1,\dots,r), ~~ {\bf e} = (1,\dots,1),
$$
$$
H'=H'({\bf v})=\left\{{\boldsymbol \eta} = 
(\eta_0, \dots, \eta_r): ~ \eta_i \in [0,1], ~ ({\boldsymbol \eta}, {\bf e}) = 1, ~ ({\boldsymbol \eta}, {\bf b})
\leq \frac{\overline{s}'-\underline{s}'}{m+1}\right\},
$$
$$
f({\bf a}) = \sum_{i=0}^r a_i \ln a_i, ~~~ {\bf a} = (a_0, a_1, \dots, a_r),
$$
$$
\rho({\bf v}) = \min_{{\boldsymbol \eta} \in H'} \{f({\boldsymbol \eta})\} - \frac{k-2}{k} f({\bf v}).
$$
Тогда 
$$
\zeta_{{\rm clique}}(k,m) \ge e^{\rho({\bf v})}.
$$}
\vskip+0.2cm
 
Следствием из теоремы является утверждение о том, что 
$$
\zeta_{{\rm clique}}(k,m) \ge \max_{r \ge 1} \max_{{\bf v} \in \Delta} e^{\rho({\bf v})}.
$$
Ясно, что мы тем самым имеем нетривиальную задачу на поиск экстремума. Мы поступим так. Сперва в параграфе 6.3 докажем теорему 10. 
Доказательство будет опираться на вероятностный метод, подобно тому, как это было с теоремами 6 и 7. Затем 
в параграфе 6.4 мы прокомментируем связь теоремы 10 с теоремами 6 и 7. В параграфе 6.5 мы обсудим задачу отыскания интересующего нас экстремума. 
Наконец, в параграфе 6.6 мы приведем таблицы результатов, которые мы получим в итоге для различных пар параметров $ k, m $. 

Заметим, что здесь мы не стали привлекать технику теории кодирования, но работаем лишь с вероятностным подходом.

\subsection{Доказательство теоремы 10}

Зафиксируем числа $ k, m $ и $ r $, а также вектор $ {\bf v} \in \Delta $. Если $ \rho({\bf v}) \le 0 $, то утверждение теоремы 10 тривиально. 
Поэтому будем считать, что $ \rho({\bf v}) > 0 $. 

Положим 
$$
{\cal V}_n = \{{\bf x} = (x_1, \dots, x_n): ~ x_i \in \{0,1, \dots, r\}, ~ 
\forall ~ j \in \{0,1, \dots, r-1\} ~~ |\{i:~x_i=j\}| = [v_jn], ~ 
$$
$$
|\{i:~x_i=r\}| = n - [v_0n] - \ldots - [v_{r-1}n]\}.
$$
Иными словами, число координат величины $ j $ в каждом векторе из $ {\cal V}_n $ асимптотически равно $ v_j n $. Нетрудно видеть, что тем самым
$$
N = |{\cal V}_n| = \left(e^{-f({\bf v})} + o(1)\right)^n.
$$

Положим, далее, 
$$
\overline{s} = \max_{{\bf x},{\bf y} \in {\cal V}_n} ({\bf x},{\bf y}), ~~~
\underline{s} = \min_{{\bf x},{\bf y} \in {\cal V}_n} ({\bf x},{\bf y}).
$$
Несложно убедиться в том, что $ \overline{s} \sim \overline{s}'n $, $ \underline{s} \sim \underline{s}'n $ (см. также \cite{Sh2}). 

Пусть $ p $ -- минимальное простое число, удовлетворяющее условию 
$ \overline{s} - (m+1)p < \underline{s} $. Тогда с учетом известных фактов из теории чисел (см. \S 2.2) имеем 
$ p \sim \frac{\overline{s}'-\underline{s}'}{m+1} $. 

Положим 
$$
{\cal A} = \left\{\sqrt{2p},\sqrt{4p}, \dots, \sqrt{2mp}\right\},
$$
$$
{\cal E}_n = \left\{\{{\bf x},{\bf y}\}:~ |{\bf x}-{\bf y}| \in {\cal A} \right\}, ~~ 
{\cal G}_n = ({\cal V}_n,{\cal E}_n) \in \mathfrak{G}_{{\cal A}}.
$$

В работе \cite{Sh1} показано, что 
$$
\alpha = \alpha({\cal G}_n) \le \sum_{{\bf u} \in H} \frac{n!}{u_0! \cdot u_1! \cdot \ldots \cdot u_r!}, 
$$
где 
$$
H=\left\{{\bf u} = (u_0, \dots, u_r): ~ u_i \in {\mathbb N}, ~ ({\bf u}, {\bf e}) = n, ~ ({\bf u}, {\bf b}) \leq p-1\right\}.
$$

Стандартный анализ с участием формулы Стирлинга дает запись 
$$ 
\alpha \le \left(e^{-\min\limits_{{\boldsymbol \eta} \in H'} f({\boldsymbol \eta})}+o(1)\right)^n. 
$$

Далее, следует применить вероятностный метод. Для этого выберем некоторые параметры. Во-первых, пусть $ \gamma \in [0,1] $ таково, что 
$ \ln \gamma < \frac{2f({\bf v})}{k} $ и что $ -\min\limits_{{\boldsymbol \eta} \in H'} \{f({\boldsymbol \eta})\} - \ln \gamma < -f({\bf v}) $. 
Последнему условию можно удовлетворить, поскольку $ \rho({\bf v}) > 0 $, т.е. 
$$
-\min\limits_{{\boldsymbol \eta} \in H'} \{f({\boldsymbol \eta})\} - \frac{2f({\bf v})}{k} < -f({\bf v}),
$$
а значит, нужно просто брать любое $ \gamma $, при котором $ \ln \gamma $ достаточно близок к величине $ \frac{2f({\bf v})}{k} $, хотя и меньше ее 
строго. Во-вторых, пусть $ c $ таково, что одновременно
$$
\ln c < -f({\bf v}), ~~~ \ln c > -\min\limits_{{\boldsymbol \eta} \in H'} \{f({\boldsymbol \eta})\} - \ln \gamma.
$$
Так, опять же, можно сделать, поскольку, ввиду неравенства $ -\min\limits_{{\boldsymbol \eta} \in H'} \{f({\boldsymbol \eta})\} - \ln \gamma < -f({\bf v}) $,
интервал 
$$
\left(-\min\limits_{{\boldsymbol \eta} \in H'} \{f({\boldsymbol \eta})\} - \ln \gamma, -f({\bf v})\right) 
$$
не пуст.     

Положим $ q = \gamma^n $, $ l = \left[c^n\right] $ (ср. \S 4.3). Построим на основе графа $ {\cal G}_n $ случайный граф с вероятностью ребра $ q $. 
Определим случайные величины $ X_l $ и $ Y_k $ в точности так же, как одноименные величины определялись в параграфе 4.3. Корректность 
и нетривиальность определения величины $ X_l $ следует из того, что 
$$ 
l = \left(e^{\ln c} + o(1)\right)^n < \left(e^{-f({\bf v})}+o(1)\right)^n = N, ~~~
l > \left(e^{-\min\limits_{{\boldsymbol \eta} \in H'} \{f({\boldsymbol \eta})\} - \ln \gamma}+o(1)\right)^n > \alpha.
$$  
Разумеется, мы всюду далее будем считать, что $ n $ достаточно велико, опуская упоминание об этом. 

Остается (см. \S 4.3) показать, что $ MX_l < \frac{1}{2} $ и $ MY_k < \frac{N}{4} $. Проводим знакомые выкладки: 
$$
MX_l \le C_N^l (1-q)^{\frac{l^2}{4\alpha}} \le e^{(c+o(1))^n-\left(\frac{\gamma c^2}{e^{-\min_{{\boldsymbol \eta} \in H'} \{f({\boldsymbol \eta})\}}}
+o(1)\right)^n}.
$$

Мы знаем, что $ \ln c > -\min\limits_{{\boldsymbol \eta} \in H'} \{f({\boldsymbol \eta})\} - \ln \gamma $. Значит, 
$ \ln c + \ln \gamma > -\min\limits_{{\boldsymbol \eta} \in H'} \{f({\boldsymbol \eta})\} $, т.е. 
$ \gamma c > e^{-\min\limits_{{\boldsymbol \eta} \in H'} \{f({\boldsymbol \eta})\}} $. Таким образом, 
$$
\frac{\gamma c^2}{e^{-\min\limits_{{\boldsymbol \eta} \in H'} \{f({\boldsymbol \eta})\}}} > c, 
$$
откуда  
$$
(c+o(1))^n-\left(\frac{\gamma c^2}{e^{-\min_{{\boldsymbol \eta} \in H'} \{f({\boldsymbol \eta})\}}}
+o(1)\right)^n \to -\infty, 
$$
а стало быть, 
$$
e^{(c+o(1))^n-\left(\frac{\gamma c^2}{e^{-\min_{{\boldsymbol \eta} \in H'} \{f({\boldsymbol \eta})\}}}
+o(1)\right)^n} < \frac{1}{2}.
$$  

В свою очередь, 
$$
MY_k \le C_N^k q^{C_k^2} \le N^k q^{C_k^2} = N \cdot N^{k-1} q^{C_k^2},
$$
и нам достаточно убедиться в том, что $ N^{k-1} q^{C_k^2} \to 0 $. В самом деле, 
$$
N^{k-1} q^{C_k^2} = \left(e^{-f({\bf v})(k-1)}+o(1)\right)^n \cdot \left(\gamma^{\frac{k(k-1)}{2}}\right)^n = 
\left(e^{-f({\bf v})(k-1)+\frac{k(k-1)}{2} \cdot \ln \gamma}+o(1)\right)^n.
$$
Но мы знаем, что $ \ln \gamma < \frac{2f({\bf v})}{k} $, а значит, 
$$
-f({\bf v})(k-1) + \frac{k(k-1)}{2} \cdot \ln \gamma < -f({\bf v})(k-1) + f({\bf v}) (k-1) = 0, 
$$
так что $ e^{-f({\bf v})(k-1)+\frac{k(k-1)}{2} \cdot \ln \gamma} \to 0 $ и $ N^{k-1} q^{C_k^2} \to 0 $.

Итак, мы показали, что 
$$
\zeta_{{\rm clique}}(k,m) \ge \frac{e^{-f({\bf v})}}{c} = e^{-f({\bf v})-\ln c}. 
$$
Беря в этом неравенстве инфимум по $ c $, для которых $ \ln c > -\min\limits_{{\boldsymbol \eta} \in H'} \{f({\boldsymbol \eta})\} - \ln \gamma $, 
и супремум по $ \gamma $, для которых $ \ln \gamma < \frac{2f({\bf v})}{k} $, получаем оценку 
$$
\zeta_{{\rm clique}}(k,m) \ge e^{-f({\bf v})+\min\limits_{{\boldsymbol \eta} \in H'} \{f({\boldsymbol \eta})\}+\frac{2f({\bf v})}{k}} = 
e^{\rho({\bf v})}. 
$$

Теорема 10 доказана.

\subsection{Небольшой комментарий к теореме 10}

Из общих соображений понятно, что теоремы 6 и 7 должны быть частными случаями теоремы 10. В действительности, так оно, конечно, и есть. 
Разница заключается в том, что в теоремах 6 и 7 мы смогли в явном виде отыскать минимум, возникающий в формулировке теоремы 10. Иными 
словами, если в теореме 10 величина $ \rho({\bf v}) $ включает в себя минимизацию по $ {\boldsymbol \eta} $, то в теоремах 6 и 7 эта 
минимизация произведена заранее. При этом в последних двух теоремах $ r = 2 $ и $ r = 3 $ соответственно. Разумеется, мы с таким же успехом 
могли в тех теоремах полагать $ r = 4 $ и т.д., рассматривая не (0,1)-векторы и не (-1,0,1)-векторы (которые с точки зрения нашей задачи 
равносильны (0,1,2)-векторам), а (0,1,2,3)-векторы и т.д. Разумеется, такое усложнение привело бы к тому, что надлежащий минимум по 
$ {\boldsymbol \eta} $ мы не смогли бы найти явно. Но суть не в том. А суть в том, что непосредственная проверка показывает отсутствие 
сколь-нибудь заметных улучшений результатов теорем 6 и 7 при переходе к большему количеству значений координат в каждом векторе 
конструкции. Именно поэтому мы выделили случай $ m = 1 $: в нем достаточно брать $ r = 2 $, $ r = 3 $, и это позволяет упростить вычисления. 

Напротив, в теореме 10 величина $ m $ произвольная, и на поверку оказывается, что чем больше $ m $, тем больше оптимальное значение $ r $. 
Очевидно, что с большими $ r $ и перебор по симплексу $ \Delta $ весьма трудоемок, а стало быть, требуются нетривиальные соображения, 
помогающие такой перебор осуществлять за разумное время. Благо для каждого $ {\bf v} \in \Delta $ нам еще нужно вычислить минимум по 
$ {\boldsymbol \eta} \in H' $, где $ H' $ довольно хитро само зависит от $ {\bf v} $. 

Аналогичную проблему в иной ситуации преодолели авторы статьи \cite{Sh2}. В следующем параграфе мы напомним экстремальную задачу, которую 
им удалось решить, и увидим, что наша задача практически с ней совпадает. В итоге нам останется лишь напомнить основные идеи из работы \cite{Sh2}, 
дабы прояснить возникновение численных данных в параграфе 6.6.

\subsection{Решение экстремальной задачи} 

Итак, нам хочется численно отыскать 
$$
\max_{r \ge 1} \max_{{\bf v} \in \Delta} e^{\min_{{\boldsymbol \eta} \in H'} \{f({\boldsymbol \eta})\} - \frac{k-2}{k} f({\bf v})}.
$$

Однако в работе \cite{Sh2} изложена технология отыскания 
$$
\max_{r \ge 1} \max_{{\bf v} \in \Delta} e^{\min_{{\boldsymbol \eta} \in H'} \{f({\boldsymbol \eta})\} - f({\bf v})}.
$$
Там слегка иные обозначения, но смысл именно такой. 

Понятно, что разница между задачами отсутствует, ведь $ k $ -- это константа. Ниже мы лишь напомним основную схему оптимизации. 

Прежде всего при каждом $ r $ ищем 
$$
\max_{{\bf v} \in \Delta} \left(\min_{{\boldsymbol \eta} \in H'} \{f({\boldsymbol \eta})\} - \frac{k-2}{k} f({\bf v})\right).
$$
Далее, замечаем, что в $ H' $ неравенство можно заменить равенством, т.е. нам нужно найти 
$$
S_{r,m,k} = \max_{{\bf v} \in \Delta} \left(\min_{{\bf s} \in \Delta, \, ({\bf s}, {\bf b})=h({\bf v})/(m+1)} f({\bf s}) - \frac{k-2}{k} f({\bf v})\right),
$$
где $ h({\bf v}) = \overline{s}'-\underline{s}' $.

Вводим векторы $ {\bf a}^i $, исходя из утверждения (предложение 2 из \cite{Sh2}): {\it функция $ h({\bf v}) $ вогнута на $ \Delta $, а неравенство 
$ h({\bf v}) \ge q $ при любом $ q $ равносильно системе из не более чем $ 2^{r-1} $ линейных неравенств $ ({\bf v},{\bf a}^i) \ge q $
с целыми неотрицательными векторами $ {\bf a}^i $, не зависящими от $ q $.} Получаем (теорема 1 из \cite{Sh2}), что 
$$
S_{r,m,k} = \max_{i=1, \dots, 2^{r-1}} \max_{0 \le p \le w(r,m)}
\left(\min_{{\bf s} \in \Delta, \, ({\bf s},{\bf b})=p} f({\bf s}) - \frac{k-2}{k} 
\min_{{\bf v} \in \Delta, \, ({\bf v},{\bf a}^i)=(m+1)p} f({\bf v})\right),
$$
где 
$$
w(r,m) = \min \left\{\frac{r-1}{2},\frac{r^2-1}{6(m+1)}\right\}.
$$

Для любого $ \lambda > 0 $ полагаем 
$$ 
{\boldsymbol \lambda} = \frac{1}{\sum\limits_{j=0}^{r-1} \lambda^j} \left(1, \lambda, \lambda^2, \dots, \lambda^{r-1}\right) \in {\mathbb R}^r.
$$
А для любого вектора $ {\bf a} = (a_1, \dots, a_r) $ с целыми неотрицательными координатами полагаем
$$
{\boldsymbol \lambda}_{{\bf a}} = \frac{1}{\sum\limits_{j=0}^{r-1} \lambda^{a_{j+1}}} 
\left(\lambda^{a_1}, \lambda^{a_2}, \dots, \lambda^{a_{r}}\right) \in {\mathbb R}^r.
$$
Оказывается, что (теорема 2 из \cite{Sh2}) 
$$
S_{r,m,k} = \max_{i=1, \dots, 2^{r-1}} \max_{0 \le p \le w(r,m)} \left(f({\boldsymbol \lambda})-\frac{k-2}{k} 
f\left({\boldsymbol \mu}_{{\bf a}^i}\right)\right),
$$
где $ \lambda $ -- единственный положительный корень полинома $ \sum\limits_{j=0}^{r-1} (j-p)z^j $, а 
$ \mu $ -- единственный положительный корень полинома $ \sum\limits_{j=0}^{r-1} \left(a_{j+1}^i-(m+1)p\right)z^{a_{j+1}^i} $. 

Такая задача решается достаточно легко. Ее численное решение описано в параграфе 4 статьи \cite{Sh2}, его мы и реализуем при получении 
таблиц из следующего параграфа, который завершает настоящий раздел, а с ним и всю нашу работу.

\subsection{Таблицы результатов}

В таблицах, которые мы приводим ниже, отражены $ m $ от 1 до 20. Для каждого $ m $ величина $ k $ меняется в стандартных пределах 1 -- 20, 25, 50, 
\dots Мы также указываем значение $ r \le 10 $, вплоть до которого приходилось искать максимумы для данных $ m $ и $ k $. Видно, что при 
$ m \le 11 $ и при $ k \to \infty $ величина искомой константы стремится ровно к той величине, которая для данного $ m $ выписана в параграфе 6.1. 
При $ m > 11 $ начинаются небольшие расхождения, связанные с тем, что мы не брали $ r > 10 $. Последнее обстоятельство обусловлено, в свою 
очередь, нехваткой мощности персонального компьютера. Если написать программу, делающую распределенные вычисления, и запустить ее на кластере, то 
можно будет брать $ r $ порядка 20, а это и есть та величина, которую использовали авторы работы \cite{Sh2}. 

Стоит отметить, что есть три клетки в первой таблице, где стоит точная единица. Это соответствует тому, что вероятностный метод дает 
экспоненциальные оценки хроматического числа не при всех $ k $. При $ m = 1 $ мы и прежде работали лишь с $ k \ge 5 $. Иными словами, 
ясно, что заведомо имеются ситуации, когда идеи теории кодирования могут помочь и здесь --- при $ m > 1 $.

\begin{longtable}{|c|c|c|c|c|c|c|c|c|c|c|}\hline
 $k \setminus m$ & $1$ & $r$ & $2$ & $r$ & $3$ & $r$ & $4$ & $r$ & $5$ & $r$ \\
\hline 
 $3$ & $1.000000$ & $0$ & $1.000000$ & $0$ & $1.001271$ & $1$ & $1.010972$ & $1$ & $1.025130$ & $2$ \\
\hline 
 $4$ & $1.000000$ & $0$ & $1.016665$ & $1$ & $1.055782$ & $2$ & $1.096186$ & $3$ & $1.134188$ & $4$ \\
\hline 
 $5$ & $1.002973$ & $1$ & $1.058718$ & $2$ & $1.122702$ & $3$ & $1.182135$ & $4$ & $1.236366$ & $5$ \\
\hline 
 $6$ & $1.018386$ & $1$ & $1.099706$ & $2$ & $1.180723$ & $3$ & $1.254357$ & $5$ & $1.321082$ & $6$ \\
\hline 
 $7$ & $1.036265$ & $1$ & $1.134785$ & $2$ & $1.229070$ & $4$ & $1.313893$ & $5$ & $1.390598$ & $7$ \\
\hline 
 $8$ & $1.052429$ & $1$ & $1.164501$ & $3$ & $1.269349$ & $4$ & $1.363230$ & $6$ & $1.448099$ & $7$ \\
\hline 
 $9$ & $1.066325$ & $1$ & $1.189781$ & $3$ & $1.303175$ & $4$ & $1.404563$ & $6$ & $1.496227$ & $7$ \\
\hline 
 $10$ & $1.078168$ & $1$ & $1.211347$ & $3$ & $1.331877$ & $4$ & $1.439591$ & $6$ & $1.536999$ & $8$ \\
\hline 
 $11$ & $1.088641$ & $2$ & $1.229900$ & $3$ & $1.356502$ & $5$ & $1.469605$ & $6$ & $1.571930$ & $8$ \\
\hline 
 $12$ & $1.098552$ & $2$ & $1.245997$ & $3$ & $1.377831$ & $5$ & $1.495581$ & $6$ & $1.602161$ & $8$ \\
\hline 
 $13$ & $1.107320$ & $2$ & $1.260076$ & $3$ & $1.396467$ & $5$ & $1.518267$ & $6$ & $1.628566$ & $8$ \\
\hline 
 $14$ & $1.115113$ & $2$ & $1.272482$ & $3$ & $1.412880$ & $5$ & $1.538243$ & $7$ & $1.651816$ & $8$ \\
\hline 
 $15$ & $1.122076$ & $2$ & $1.283491$ & $3$ & $1.427439$ & $5$ & $1.555959$ & $7$ & $1.672439$ & $8$ \\
\hline 
 $16$ & $1.128325$ & $2$ & $1.293319$ & $3$ & $1.440437$ & $5$ & $1.571776$ & $7$ & $1.690852$ & $8$ \\
\hline 
 $17$ & $1.133962$ & $2$ & $1.302145$ & $3$ & $1.452111$ & $5$ & $1.585980$ & $7$ & $1.707389$ & $9$ \\
\hline 
 $18$ & $1.139068$ & $2$ & $1.310113$ & $3$ & $1.462651$ & $5$ & $1.598803$ & $7$ & $1.722321$ & $9$ \\
\hline 
 $19$ & $1.143714$ & $2$ & $1.317370$ & $4$ & $1.472213$ & $5$ & $1.610438$ & $7$ & $1.735870$ & $9$ \\
\hline 
 $20$ & $1.147957$ & $2$ & $1.324001$ & $4$ & $1.480927$ & $5$ & $1.621040$ & $7$ & $1.748218$ & $9$ \\
\hline 
 $25$ & $1.164612$ & $2$ & $1.349947$ & $4$ & $1.514987$ & $5$ & $1.662487$ & $7$ & $1.796502$ & $9$ \\
\hline 
 $50$ & $1.200435$ & $2$ & $1.405522$ & $4$ & $1.587854$ & $6$ & $1.751161$ & $8$ & $1.899879$ & $9$ \\
\hline 
 $100$ & $1.219590$ & $2$ & $1.435237$ & $4$ & $1.626797$ & $6$ & $1.798568$ & $8$ & $1.955191$ & $10$ \\
\hline 
 $1000$ & $1.237532$ & $2$ & $1.463131$ & $4$ & $1.663355$ & $6$ & $1.843092$ & $8$ & $2.007171$ & $10$ \\
\hline 
 $10000$ & $1.239363$ & $2$ & $1.465982$ & $4$ & $1.667092$ & $6$ & $1.847644$ & $8$ & $2.012487$ & $10$ \\
\hline 
 $100000$ & $1.239546$ & $2$ & $1.466267$ & $4$ & $1.667467$ & $6$ & $1.848100$ & $8$ & $2.013020$ & $10$ \\
\hline 
\end{longtable}

\begin{longtable}{|c|c|c|c|c|c|c|c|c|c|c|}\hline
 $k\setminus m$ & $6$ & $r$ & $7$ & $r$ & $8$ & $r$ & $9$ & $r$ & $10$ & $r$ \\
\hline 
 $3$ & $1.040369$ & $3$ & $1.055470$ & $3$ & $1.070035$ & $4$ & $1.083948$ & $5$ & $1.097196$ & $5$ \\
\hline 
 $4$ & $1.169410$ & $5$ & $1.202065$ & $6$ & $1.232460$ & $7$ & $1.260886$ & $8$ & $1.287592$ & $9$ \\
\hline 
 $5$ & $1.286031$ & $6$ & $1.331837$ & $8$ & $1.374384$ & $9$ & $1.414155$ & $10$ & $1.451537$ & $10$ \\
\hline 
 $6$ & $1.382065$ & $7$ & $1.438305$ & $9$ & $1.490586$ & $10$ & $1.539514$ & $10$ & $1.585569$ & $10$ \\
\hline 
 $7$ & $1.460707$ & $8$ & $1.525420$ & $9$ & $1.585647$ & $10$ & $1.642088$ & $10$ & $1.695284$ & $10$ \\
\hline 
 $8$ & $1.525716$ & $8$ & $1.597430$ & $10$ & $1.664250$ & $10$ & $1.726939$ & $10$ & $1.786094$ & $10$ \\
\hline 
 $9$ & $1.580123$ & $9$ & $1.657711$ & $10$ & $1.730074$ & $10$ & $1.798032$ & $10$ & $1.862219$ & $10$ \\
\hline 
 $10$ & $1.626218$ & $9$ & $1.708798$ & $10$ & $1.785883$ & $10$ & $1.858336$ & $10$ & $1.926821$ & $10$ \\
\hline 
 $11$ & $1.665716$ & $9$ & $1.752588$ & $10$ & $1.833738$ & $10$ & $1.910065$ & $10$ & $1.982261$ & $10$ \\
\hline 
 $12$ & $1.699908$ & $10$ & $1.790506$ & $10$ & $1.875191$ & $10$ & $1.954891$ & $10$ & $2.030319$ & $10$ \\
\hline 
 $13$ & $1.729777$ & $10$ & $1.823641$ & $10$ & $1.911426$ & $10$ & $1.994086$ & $10$ & $2.072354$ & $10$ \\
\hline 
 $14$ & $1.756084$ & $10$ & $1.852831$ & $10$ & $1.943356$ & $10$ & $2.028635$ & $10$ & $2.109415$ & $10$ \\
\hline 
 $15$ & $1.779423$ & $10$ & $1.878734$ & $10$ & $1.971698$ & $10$ & $2.059309$ & $10$ & $2.142327$ & $10$ \\
\hline 
 $16$ & $1.800265$ & $10$ & $1.901871$ & $10$ & $1.997019$ & $10$ & $2.086718$ & $10$ & $2.171743$ & $10$ \\
\hline 
 $17$ & $1.818987$ & $10$ & $1.922658$ & $10$ & $2.019773$ & $10$ & $2.111354$ & $10$ & $2.198188$ & $10$ \\
\hline 
 $18$ & $1.835895$ & $10$ & $1.941435$ & $10$ & $2.040329$ & $10$ & $2.133615$ & $10$ & $2.222086$ & $10$ \\
\hline 
 $19$ & $1.851238$ & $10$ & $1.958477$ & $10$ & $2.058990$ & $10$ & $2.153826$ & $10$ & $2.243787$ & $10$ \\
\hline 
 $20$ & $1.865223$ & $10$ & $1.974012$ & $10$ & $2.076004$ & $10$ & $2.172256$ & $10$ & $2.263578$ & $10$ \\
\hline 
 $25$ & $1.919929$ & $10$ & $2.034805$ & $10$ & $2.142603$ & $10$ & $2.244423$ & $10$ & $2.341099$ & $10$ \\
\hline 
 $50$ & $2.037155$ & $10$ & $2.165191$ & $10$ & $2.285568$ & $10$ & $2.399459$ & $10$ & $2.507751$ & $10$ \\
\hline 
 $100$ & $2.099938$ & $10$ & $2.235086$ & $10$ & $2.362273$ & $10$ & $2.482705$ & $10$ & $2.597292$ & $10$ \\
\hline 
 $1000$ & $2.158974$ & $10$ & $2.300852$ & $10$ & $2.434486$ & $10$ & $2.561114$ & $10$ & $2.681663$ & $10$ \\
\hline 
 $10000$ & $2.165015$ & $10$ & $2.307583$ & $10$ & $2.441879$ & $10$ & $2.569143$ & $10$ & $2.690304$ & $10$ \\
\hline 
 $100000$ & $2.165620$ & $10$ & $2.308257$ & $10$ & $2.442620$ & $10$ & $2.569948$ & $10$ & $2.691170$ & $10$ \\
\hline 
\end{longtable}

\begin{longtable}{|c|c|c|c|c|c|c|c|c|c|c|}\hline
 $k\setminus m$ & $11$ & $r$ & $12$ & $r$ & $13$ & $r$ & $14$ & $r$ & $15$ & $r$ \\
\hline 
 $3$ & $1.109804$ & $6$ & $1.121815$ & $7$ & $1.133272$ & $7$ & $1.144220$ & $10$ & $1.154699$ & $10$ \\
\hline 
 $4$ & $1.312789$ & $10$ & $1.336654$ & $10$ & $1.359334$ & $10$ & $1.380954$ & $10$ & $1.401621$ & $10$ \\
\hline 
 $5$ & $1.486840$ & $10$ & $1.520320$ & $10$ & $1.552185$ & $10$ & $1.582609$ & $10$ & $1.611740$ & $10$ \\
\hline 
 $6$ & $1.629130$ & $10$ & $1.670504$ & $10$ & $1.709943$ & $10$ & $1.747657$ & $10$ & $1.783821$ & $10$ \\
\hline 
 $7$ & $1.745667$ & $10$ & $1.793584$ & $10$ & $1.839319$ & $10$ & $1.883104$ & $10$ & $1.925138$ & $10$ \\
\hline 
 $8$ & $1.842182$ & $10$ & $1.895581$ & $10$ & $1.946598$ & $10$ & $1.995487$ & $10$ & $2.042459$ & $10$ \\
\hline 
 $9$ & $1.923133$ & $10$ & $1.981176$ & $10$ & $2.036674$ & $10$ & $2.089894$ & $10$ & $2.141061$ & $10$ \\
\hline 
 $10$ & $1.991863$ & $10$ & $2.053882$ & $10$ & $2.113218$ & $10$ & $2.170152$ & $10$ & $2.224914$ & $10$ \\
\hline 
 $11$ & $2.050870$ & $10$ & $2.116326$ & $10$ & $2.178982$ & $10$ & $2.239126$ & $10$ & $2.296999$ & $10$ \\
\hline 
 $12$ & $2.102037$ & $10$ & $2.170490$ & $10$ & $2.236042$ & $10$ & $2.298989$ & $10$ & $2.359574$ & $10$ \\
\hline 
 $13$ & $2.146803$ & $10$ & $2.217892$ & $10$ & $2.285990$ & $10$ & $2.351400$ & $10$ & $2.414370$ & $10$ \\
\hline 
 $14$ & $2.186283$ & $10$ & $2.259706$ & $10$ & $2.330059$ & $10$ & $2.397649$ & $10$ & $2.462731$ & $10$ \\
\hline 
 $15$ & $2.221351$ & $10$ & $2.296853$ & $10$ & $2.369216$ & $10$ & $2.438750$ & $10$ & $2.505712$ & $10$ \\
\hline 
 $16$ & $2.252699$ & $10$ & $2.330066$ & $10$ & $2.404231$ & $10$ & $2.475507$ & $10$ & $2.544154$ & $10$ \\
\hline 
 $17$ & $2.280886$ & $10$ & $2.359933$ & $10$ & $2.435722$ & $10$ & $2.508568$ & $10$ & $2.578734$ & $10$ \\
\hline 
 $18$ & $2.306362$ & $10$ & $2.386932$ & $10$ & $2.464192$ & $10$ & $2.538460$ & $10$ & $2.610001$ & $10$ \\
\hline 
 $19$ & $2.329498$ & $10$ & $2.411454$ & $10$ & $2.490052$ & $10$ & $2.565614$ & $10$ & $2.638405$ & $10$ \\
\hline 
 $20$ & $2.350601$ & $10$ & $2.433823$ & $10$ & $2.513644$ & $10$ & $2.590387$ & $10$ & $2.664320$ & $10$ \\
\hline 
 $25$ & $2.433281$ & $10$ & $2.521481$ & $10$ & $2.606109$ & $10$ & $2.687495$ & $10$ & $2.765914$ & $10$ \\
\hline 
 $50$ & $2.611128$ & $10$ & $2.710124$ & $10$ & $2.805166$ & $10$ & $2.896599$ & $10$ & $2.984709$ & $10$ \\
\hline 
 $100$ & $2.706733$ & $10$ & $2.811571$ & $10$ & $2.912242$ & $10$ & $3.009098$ & $10$ & $3.102433$ & $10$ \\
\hline 
 $1000$ & $2.796842$ & $10$ & $2.907207$ & $10$ & $3.013196$ & $10$ & $3.115171$ & $10$ & $3.213433$ & $10$ \\
\hline 
 $10000$ & $2.806072$ & $10$ & $2.917003$ & $10$ & $3.023538$ & $10$ & $3.126038$ & $10$ & $3.224805$ & $10$ \\
\hline 
 $100000$ & $2.806998$ & $10$ & $2.917985$ & $10$ & $3.024575$ & $10$ & $3.127127$ & $10$ & $3.225945$ & $10$ \\
\hline 
\end{longtable}

\begin{longtable}{|c|c|c|c|c|c|c|c|c|c|c|}\hline
 $k\setminus m$ & $16$ & $r$ & $17$ & $r$ & $18$ & $r$ & $19$ & $r$ & $20$ & $r$ \\
\hline 
 $3$ & $1.164748$ & $10$ & $1.174400$ & $9$ & $1.183687$ & $10$ & $1.192635$ & $10$ & $1.201269$ & $10$ \\
\hline 
 $4$ & $1.421424$ & $10$ & $1.440444$ & $10$ & $1.458746$ & $10$ & $1.476392$ & $10$ & $1.493431$ & $10$ \\
\hline 
 $5$ & $1.639701$ & $10$ & $1.666600$ & $10$ & $1.692528$ & $10$ & $1.717566$ & $10$ & $1.741784$ & $10$ \\
\hline 
 $6$ & $1.818582$ & $10$ & $1.852068$ & $10$ & $1.884388$ & $10$ & $1.915635$ & $10$ & $1.945892$ & $10$ \\
\hline 
 $7$ & $1.965584$ & $10$ & $2.004585$ & $10$ & $2.042260$ & $10$ & $2.078714$ & $10$ & $2.114038$ & $10$ \\
\hline 
 $8$ & $2.087693$ & $10$ & $2.131339$ & $10$ & $2.173527$ & $10$ & $2.214368$ & $10$ & $2.253957$ & $10$ \\
\hline 
 $9$ & $2.190361$ & $10$ & $2.237954$ & $10$ & $2.283974$ & $10$ & $2.328537$ & $10$ & $2.371745$ & $10$ \\
\hline 
 $10$ & $2.277700$ & $10$ & $2.328674$ & $10$ & $2.377976$ & $10$ & $2.425726$ & $10$ & $2.472028$ & $10$ \\
\hline 
 $11$ & $2.352800$ & $10$ & $2.406697$ & $10$ & $2.458835$ & $10$ & $2.509336$ & $10$ & $2.558307$ & $10$ \\
\hline 
 $12$ & $2.418003$ & $10$ & $2.474448$ & $10$ & $2.529056$ & $10$ & $2.581952$ & $10$ & $2.633247$ & $10$ \\
\hline 
 $13$ & $2.475109$ & $10$ & $2.533792$ & $10$ & $2.590568$ & $10$ & $2.645566$ & $10$ & $2.698897$ & $10$ \\
\hline 
 $14$ & $2.525513$ & $10$ & $2.586176$ & $10$ & $2.644869$ & $10$ & $2.701724$ & $10$ & $2.756855$ & $10$ \\
\hline 
 $15$ & $2.570315$ & $10$ & $2.632739$ & $10$ & $2.693139$ & $10$ & $2.751646$ & $10$ & $2.808377$ & $10$ \\
\hline 
 $16$ & $2.610388$ & $10$ & $2.674391$ & $10$ & $2.736318$ & $10$ & $2.796303$ & $10$ & $2.854465$ & $10$ \\
\hline 
 $17$ & $2.646437$ & $10$ & $2.711861$ & $10$ & $2.775162$ & $10$ & $2.836478$ & $10$ & $2.895928$ & $10$ \\
\hline 
 $18$ & $2.679033$ & $10$ & $2.745743$ & $10$ & $2.810288$ & $10$ & $2.872807$ & $10$ & $2.933421$ & $10$ \\
\hline 
 $19$ & $2.708646$ & $10$ & $2.776525$ & $10$ & $2.842201$ & $10$ & $2.905813$ & $10$ & $2.967485$ & $10$ \\
\hline 
 $20$ & $2.735666$ & $10$ & $2.804611$ & $10$ & $2.871318$ & $10$ & $2.935928$ & $10$ & $2.998564$ & $10$ \\
\hline 
 $25$ & $2.841592$ & $10$ & $2.914723$ & $10$ & $2.985476$ & $10$ & $3.053995$ & $10$ & $3.120412$ & $10$ \\
\hline 
 $50$ & $3.069738$ & $10$ & $3.151892$ & $10$ & $3.231353$ & $10$ & $3.308281$ & $10$ & $3.382821$ & $10$ \\
\hline 
 $100$ & $3.192495$ & $10$ & $3.279500$ & $10$ & $3.363639$ & $10$ & $3.445081$ & $10$ & $3.523978$ & $10$ \\
\hline 
 $1000$ & $3.308239$ & $10$ & $3.399814$ & $10$ & $3.488357$ & $10$ & $3.574045$ & $10$ & $3.657042$ & $10$ \\
\hline 
 $10000$ & $3.320097$ & $10$ & $3.412139$ & $10$ & $3.501132$ & $10$ & $3.587256$ & $10$ & $3.670672$ & $10$ \\
\hline 
 $100000$ & $3.321285$ & $10$ & $3.413375$ & $10$ & $3.502413$ & $10$ & $3.588580$ & $10$ & $3.672038$ & $10$ \\
\hline 
\end{longtable}

\maketitle

%\begin{large}

%\end{large}

\end{document}